\documentclass[a4paper,11pt]{amsart}
\usepackage{amssymb,amsfonts,amsmath}

\def\d{\delta}

\def\C{\mathbb{C}}
\def\c2{\mathbb{C}^2}
\def\R{\mathbb{R}}

\def\Z{\mathbb{Z}}
\def\N{\mathbb{N}}
\def\P{\mathbb{P}}

\def\1{\bold{1}}
\def\B{\mathbb{B}}

\def\a{\alpha}
\def\b{\beta}
\def\e{\varepsilon}
\def\l{\lambda}
\def\f{\varphi}
\def\g{\gamma}

\def\p{\psi}
\def\r{\varrho}
\def\o{\omega}

\def\D{\overline{\partial}}

\newtheorem{lem}{Lemma}[section]
\newtheorem{pro}[lem]{Proposition}
\newtheorem{defi}[lem]{Definition}
\newtheorem{def/not}[lem]{Definition/Notations}

\newtheorem{thm}[lem]{Theorem}
\newtheorem{cor}[lem]{Corollary}
\newtheorem{rqe}[lem]{Remark}

\newtheorem{exa}[lem]{Example}

\newenvironment{proof3.1}
{\noindent {\it{Proof of theorem 3.1}}}{$\Box$ \linebreak[4]}

\newenvironment{warning}
{\vskip .2cm \noindent {{\bf Warning.}}}{\hfill\break}

\begin{document}

\title[Intrinsic capacities on compact K\"ahler manifolds]
{Intrinsic capacities on compact K\"ahler manifolds}

\author{Vincent GUEDJ \& Ahmed ZERIAHI}

\begin{abstract}
We study fine properties of quasiplurisubharmonic functions on compact K\"ahler
manifolds. We define and study several intrinsic capacities which characterize
pluripolar sets and show that locally pluripolar sets are globally
"quasi-pluripolar". 
\end{abstract}

\maketitle


\section*{Introduction}

Since the fundamental work of Bedford and Taylor [4],[5],
several authors have
developed a "Pluripotential theory" in domains of $\C^n$
(or of Stein manifolds). This theory is devoted to the fine study
of plurisubharmonic (psh) functions and can be seen as a non-linear
generalization of the classical potential theory (in one complex variable),
where subharmonic functions and the Laplace operator $\Delta$ are replaced
by psh functions and the complex Monge-Amp\`ere operator $(dd^c)^n$.
Here $d,d^c$ denote the real differential operators
$d:=\partial+\overline{\partial}$,
$d^c:=\frac{i}{2 \pi}[\overline{\partial}-\partial]$ so that
$dd^c=\frac{i}{\pi} \partial \overline{\partial}$ ; the normalization
being chosen so that the positive measure 
$\left( dd^c \frac{1}{2}\log[1+||z||^2] \right)^n$ has total
mass $1$ in $\C^n$.
We refer the reader to [3],[9],[28], [29] for a
survey of this local theory.

Our aim here is to develop a global Pluripotential theory in the
context of compact K\"ahler manifolds. It follows from the 
maximum principle that there are no psh functions (except
constants) on a compact complex manifold $X$. However there are 
usually plenty of positive closed currents of bidegree $(1,1)$
(we refer the reader to [17], chapter 3, for basic facts
on positive currents). Given $\o$ a real closed smooth form
of bidegree $(1,1)$ on $X$, we may consider
every positive closed current $\o'$ of bidegree $(1,1)$ on $X$
which is cohomologous to $\o$. When $X$ is K\"ahler, it follows
from the "$dd^c$-lemma" 
that $\o'$ can be written as $\o'=\o+dd^c \f$, where $\f$ is
a function which is integrable with respect to any smooth volume form
on $X$. Such a function $\f$ 
will be called $\o$-plurisubharmonic ($\o$-psh for short).
It is globally defined on $X$ and locally given as the sum
of a psh and a smooth function. 
We let $PSH(X,\o)$ denote the set of $\o$-psh functions.
Such functions were introduced by Demailly,
who call them {\it quasiplurisubharmonic} (qpsh). These are the main
objects of study in this article.

There are several motivations to study qpsh functions on compact 
K\"ahler manifolds. First of all they arise naturally in complex analytic
geometry as positive singular metrics of holomorphic line bundles
(see section 4) whose study is central to several 
questions of complex algebraic geometry.
Solving Monge-Amp\`ere equations associated to $\o$-psh functions 
has been used to produce metrics with prescribed singularities (see [15]). 
It is also related to the
existence of canonical metrics in K\"ahler geometry (see [39]).
Important contributions have been made in this direction by Blocki (see [6], [7]) and also by Kolodziej using techniques from local Pluripotential theory (see [30], [31]). 
Quasiplurisubharmonic functions have also been used in [23] to define a notion
of $\o$-polynomial convexity and study the fine approximation of 
positive currents by rational divisors. Last but not least, such functions are
of constant use in complex dynamics in several variables (see [19], [20],
[24], [25], [35]).

It seems to us appropriate to develop a theory of qpsh functions
of its own rather than view these functions as particular cases of
the local theory. 
Although the two theories look quite similar, there are important differences
which make the "compact theory" both simpler and more difficult than the local
one. Here are some examples:
\begin{itemize}
\item  There is no pluriharmonic functions (except constants)
on a compact manifold, hence each $\o$-psh function $\f$ is canonically
associated (up to normalization) to its curvature current
$\o_{\f}:=\o+dd^c \f \geq 0$. This yields compactness
properties of subsets of $PSH(X,\o)$ 
(see section 1) which are quite useful
(e.g. in complex dynamics, see section 6.2). 

\item Integration by parts (of constant use in such theories) is quite simple
  in the compact setting since there is no boundary.
As an illustration, we obtain transparent proofs of Chern-Levine-Nirenberg
type inequalities (see example 1.8 and section 2).
A successful application of this simple observation has been made in complex
dynamics in [25].

\item On the other hand one looses homogeneity of Monge-Amp\`ere operators in
  the compact setting. They do have uniformly bounded mass (by Stokes
  theorem), but there is no performing "comparison principle", which is a key
  tool in the local theory. This is a source of difficulty when, for example,
  one wishes to solve Monge-Amp\`ere equations on compact manifolds (see [30],
  [26]). 
\end{itemize}
We shall develop our study in a series of articles. 
In the present one we define and study several intrinsic capacities which we
shall use in our forthcoming articles.
\vskip.1cm

Let us now describe more precisely the contents of the article.

\noindent In {\it section 1} we define $\o$-psh functions and gather useful
facts about them (especially compacity results such as proposition 1.7). For
locally bounded $\o$-psh functions $\f$ we define the complex Monge-Amp\`ere
operator $\o_{\f}^n$ in {\it section 2}. We establish Chern-Levine-Nirenberg
inequalities (proposition 2.1) and study the "Monge-Amp\`ere capacity"
$Cap_{\o}$ (definition 2.4). As in the local theory, $\o$-psh functions are
quasicontinuous with respect to $Cap_{\o}$ (corollary 2.8). The capacity
$Cap_{\o}$ is comparable to the local Monge-Amp\`ere capacity of Bedford and
Taylor (proposition 2.10) and moreover enjoys invariance properties
(proposition 2.5). 
In {\it section 3} we define a relative extremal function 
$h_{E,\o}^*$ and establish a useful formula (theorem 3.2)
$$
Cap_{\o}^*(E)=\int_X (-h_{E,\o}^*) \left( \o+dd^c h_{E,\o}^*\right)^n.
$$
This is the global version of the fundamental local formula
of Bedford-Taylor [5],
$Cap(E,\omega)=\int_{\Omega} (dd^c u_E^*)^n$.

In {\it section 4} we study yet another capacity (the
Alexander capacity $T_{\o}$, definition 4.7) which is defined by means of a
(global) 
extremal function (definition 4.1). When $\o$ is a Hodge form, it can be
defined as well in terms of Tchebychev constants: these are the contents of 
{\it section 5} (theorem 5.2) where we further give a geometrical
interpretation of $T_{\o}$ when $X=\C\P^n$ is the complex projective space and
$\o$ is the Fubini-Study K\"ahler form (theorem 5.4), following Alexander's
work [1]. In {\it section 6} we
show that locally pluripolar sets can be defined by $\o$-psh functions when
$\o$ is K\"ahler: this is our version of a result of Josefson
(theorem 6.2). We then give an application in complex dynamics which
illustrates how invariance properties of these capacities can be used.
Finally in an Appendix we show how to globally regularize $\o$-psh functions,
following ideas of Demailly.

\vskip.1cm

This paper lies at the border of Complex Analysis and Complex Geometry.
We have tried to make it accessible to mathematicians from both sides.
This has of course some consequences for the style of presentation.
We have included proofs of some results which may be seen as
consequences of results from the local pluripotential theory.
We have spent some efforts defining, regularizing and approximating
positive singular metrics of holomorphic line bundles, although
some of these facts may be considered as classical by complex geometers.
Altogether we hope the paper is essentially self-contained. Our efforts
will not be vain if for instance we have convinced specialists of the (local)
pluripotential theory that the right point of view in studying the Lelong 
class ${\mathcal L}(\C^n)$ of psh functions with logarithmic growth
in $\C^n$ is to consider qpsh functions on the complex projective
space $\C\P^n$. We also think this paper should be useful to people working in
complex dynamics in several variables where pluripotential theory has become
an important tool.

\begin{warning}
In the whole paper {\it positivity} (like e.g. in {\it positive} metric
and {\it positive} current) has to be understood in the weak
(french, i.e. non-negativity) sense of currents, except when we
talk of a {\it positive} line bundle $L$, in which case it means that
$L$ admits a smooth metric whose curvature is a K\"ahler form.
\end{warning}


\section{Quasiplurisubharmonic functions}

In the sequel, unless otherwise specified, $L^p$-norms will always be computed
with respect to a fixed volume form on $X$, which is a 
{\it compact connected K\"ahler} manifold.
Let $\o$ be a closed real current of bidegree $(1,1)$ on $X$. We assume
throughout the article that $\o$ has continuous local potentials.

\begin{defi} Set
$$
PSH(X,\o):=\{ \f \in L^1(X,\R \cup \{ -\infty \}) \, / \, 
dd^c \f \geq -\o \text{ and } \f
\text{ is } u.s.c. \}.$$
The set $PSH(X,\o)$ is the set of "$\o$-plurisubharmonic" functions.
\end{defi}

Observe that $PSH(X,\o)$ is non empty if and only if there exists
a positive closed current of bidegree $(1,1)$ on $X$ which is
cohomologous to $\o$. One then says that the cohomology class
$[\o]$ is {\it pseudoeffective}. In the sequel we always assume
this property holds.
We also always assume that $\o$ has continuous local potentials.
This guarantees that $\o$-plurisubharmonic functions ($\o$-psh for short)
are upper semi-continuous (u.s.c.), so they are
locally hence globally bounded from above.
We endow $PSH(X,\o)$ with the $L^1$-topology. Observe that $PSH(X,\o)$
is a closed subspace of $L^1(X)$.

\begin{exa}
The most fundamental example which may serve as a guideline to everything
that follows is the case where $X=\C\P^n$ is the complex projective space
and $\o=\o_{FS}$ is the Fubini-Study K\"ahler form. There is then a 
1-to-1 correspondence between $PSH(\C\P^n,\o_{FS})$ and the Lelong class
$$
{\mathcal L}(\C^n):=\left\{ \p \in PSH(\C^n) \, / \, \p(z) \leq \frac{1}{2}
  \log [1+|z|^2]+C_{\p} \right\}
$$
which is given by the natural mapping
$$
\p \in {\mathcal L}(\C^n) \mapsto \f(x)=\left\{ 
\begin{array}{ccc} 
\p(x)-\frac{1}{2} \log [1+|x|^2] & \text{ if } & x \in \C^n \\
 \overline{\lim}_{y \in \C^n \rightarrow x} (\p(y)-\frac{1}{2} \log [1+|y|^2])
 & \text{ if } & x \in H_{\infty} ,
\end{array} \right.
$$
where $H_{\infty}$ denotes the hyperplane at infinity. One can easily show
that this mapping is bicontinuous for the $L_{loc}^1$ topology.
\end{exa}

The Lelong class ${\mathcal L}(\C^n)$ of plurisubharmonic functions
with logarithmic growth in $\C^n$ has been intensively studied 
in the last thirty years.
It seems to us that the properties of ${\mathcal L}(\C^n)$ are more easily
seen when ${\mathcal L}(\C^n)$ is viewed as $PSH(\C\P^n,\o_{FS})$.
Further we shall see hereafter that the class $PSH(X,\o)$ of 
$\o$-psh functions enjoys several properties of ${\mathcal L}(\C^n)$
when $\o$ is K\"ahler. We start by observing (proposition 1.3.1 \& 1.3.2
below) that $PSH(X,\o)$ and $PSH(X,\o')$ are comparable if $\o,\o'$
are both K\"ahler.

\begin{pro} \text{ }

1) If $\o_1 \leq \o_2$ then $PSH(X,\o_1) \subset PSH(X,\o_2)$.

2) $\forall A \in \R_+^*$, $PSH(X,A \o)= A \cdot PSH(X,\o)$.

3) If $\o'$ is cohomologous to $\o$, $\o'=\o+dd^c \chi$, then
$$
PSH(X,\o')=PSH(X,\o)+\chi.
$$

4) If $\f,\p \in PSH(X,\o)$ then 
$$
\max(\f,\p) \; , \; \frac{\f+\p}{2} \; , \; \log[e^{\f}+e^{\p}] \in PSH(X,\o)
$$
\end{pro}

\begin{proof}
Assertions 1),2),3) follow straightforwardly from the definition.
Observe that 1.3.4 says that $PSH(X,\o)$
is a convex set which is stable under taking maximum and also under
the operation $(\f,\p) \mapsto \log[e^{\f}+e^{\p}]$. These are
all consequences of the corresponding local properties of psh functions.
We nevertheless give a proof, in the spirit of this article. 
That $(\f+\p)/2 \in PSH(X,\o)$ follows by linearity.
The latter assertion is a consequence of the
following computation
$$
dd^c \log[e^{\f}+e^{\p}] = \frac{e^{\f}dd^c \f +e^{\p} dd^c \p}{e^{\f}+e^{\p}}
+\frac{e^{\f+\p} d(\f-\p) \wedge d^c (\f-\p)}{[e^{\f}+e^{\p}]^2},
$$
using that $df \wedge d^c f \geq 0$. This computation makes sense
if for instance $\f,\p$ are smooth. The general case follows then by
regularizing $\f,\p$ (see Appendix). Finally observe that
$\max(\f,\p)=\lim j^{-1}\log [e^{j\f}+e^{j\p}] \in PSH(X,\o)$.
\end{proof}

It follows from 1.3.3 that $PSH(X,\o)$ essentially depends on the 
cohomology class $[\o]$. In the same vein we have the following:

\begin{pro}
Let ${\mathcal T}_{[\o]}(X)$ denote the set of positive closed currents $\o'$
of bidegree $(1,1)$ on $X$ which are cohomologous to $\o$. Then
$$
PSH(X,\o) \simeq {\mathcal T}_{[\o]}(X) \oplus \R.
$$
\end{pro}

\begin{proof}
The mapping 
$$
\Phi: \f \in PSH(X,\o) \mapsto \o_{\f}:=\o+dd^c \f \in {\mathcal T}_{[\o]}(X)
$$
is a continuous affine mapping whose kernel consists of constants mappings: 
indeed $\o_{\f}=\o_{\p}$ implies
that $\f-\p$ is pluriharmonic hence constant by the maximum principle.
Moreover $\Phi$ is surjective: if $\o' \geq 0$ is cohomologous to 
$\o$ then $\o'=\o+dd^c \f$ for some $\f \in L^1(X,\R)$ -this is the
celebrated $dd^c$-lemma on K\"ahler manifolds (see e.g. lemma 8.6,
chapter VI in [17]).
Thus $\f$ coincides almost everywhere with a function of $PSH(X ,\o)$
and $\o'=\Phi(\f)$.
\end{proof}

\begin{rqe}
The size of $PSH(X,\o)$ is therefore related to that of 
${\mathcal T}_{[\o]}(X)$ hence only depends on the positivity of the
cohomology class $[\o]$. The more positive $[\o]$, the bigger
$PSH(X,\o)$. 

When $[\o]$ is K\"ahler then $PSH(X,\o)$ is large: if e.g.
$\chi$ is any ${\mathcal C}^2$-function on $X$ then 
$\e \chi \in PSH(X,\o)$ for $\e>0$ small enough.
We will see (theorem 6.2) that $PSH(X,\o)$ characterizes
locally pluripolar sets when $[\o]$ is K\"ahler. It follows from
proposition 1.3 that $PSH(X,\o)$ and $PSH(X,\o')$ have the same "size"
if $\o$ and $\o'$ are both K\"ahler.

Note on the other
hand that $PSH(X,\o) \simeq \R$ when $\o$ is cohomologous to
$[E]$, the current of integration
along the exceptional divisor of a smooth blow up. Indeed let
$\pi:X \rightarrow \tilde{X}$ be a blow up with smooth center $Y$,
$ codim_{\C} Y \geq 2$ (see e.g. chapter 2 of [17]  
for the definition of 
blow-ups). Let $E=\pi^{-1}(Y)$ denote the exceptional divisor
and $\o=[E]$ be the current of integration along $E$.
If $\f \in PSH(X,[E])$ then $dd^c (\f \circ \pi^{-1}) \geq 0$
in $\tilde{X} \setminus Y$. Since $ codim_{\C} Y \geq 2$,
$\f \circ \pi^{-1}$ extends trivially through $Y$ has a
global psh function on $\tilde{X}$. By the maximum principle
$\f \circ \pi^{-1}$ is constant hence so is $\f$.
Alternatively there is no positive closed current of bidegree
$(1,1)$ on $X$ which is cohomologous to $[E]$ except $[E]$
itself.
\end{rqe}

It follows from previous proposition that any set of "normalized" 
$\o$-psh functions is in 1-to-1 correspondence with 
${\mathcal T}_{[\o]}(X)$ which is compact for the weak topology of currents.
This is the key to several results to follow: normalized $\o$-psh
functions form a compact family in $L^1(X)$.

\begin{pro}
Let $(\f_j) \in PSH(X,\o)^{\N}$. 

1) If $(\f_j)$ is uniformly bounded from above on $X$, then either
$\f_j$ converges uniformly to $-\infty$ on $X$ or the sequence 
$(\f_j)$ is relatively compact in $L^1(X)$.

2) If $\f_j \rightarrow \f$ in $L^1(X)$, then $\f$ coincides almost everywhere
   with a unique function $\f^* \in PSH(X,\o)$. Moreover
$$
\sup_X \f^*=\lim_{j \rightarrow +\infty} \sup_X \f_j.
$$

3) In particular if $\f_j$ is decreasing, then either 
$\f_j \rightarrow -\infty$ or $\f=\lim \f_j \in PSH(X,\o)$.
Similarly, if $\f_j$ is increasing and uniformly bounded from above then
$\f:=(\lim \f_j)^* \in PSH(X,\o)$, where ${\cdot}^*$ denotes the
upper-semi-continuous regularization.
\end{pro}

\begin{proof}
This is a straighforward consequence of the analogous local result for
sequences of psh functions. We refer the reader to 
[17], chapter 1, for a proof. Note that 1.6.2 is a special case
of a celebrated lemma attributed to Hartogs.
\end{proof}

The next result is quite useful (see [42], [43] for a systematic use).

\begin{pro}
The family 
$$
{\mathcal F}_0:=\{ \f \in PSH(X,\o) \, / \, \sup_X \f=0 \}
$$
is a compact subset of $PSH(X,\o)$.

If $\mu$ is a probability measure such that 
$PSH(X,\o) \subset L^1(\mu)$ then
$$
{\mathcal F}_{\mu}:=\{ \f \in PSH(X,\o) \, / \, \int_X \f d\mu=0 \}
$$
is a relatively compact subset of $PSH(X,\o)$. 
In particular there exists $C_{\mu}$
such that $\forall \f \in PSH(X,\o)$,
$$
-C_{\mu}+ \sup_X \f \leq \int_X \f d\mu \leq \sup_X \f .
$$
\end{pro}

\begin{proof}
It follows straightforwardly from proposition 1.6.1 that 
${\mathcal F}_0$ is a relatively compact subset of $PSH(X,\o)$.
Moreover ${\mathcal F}_0$ is closed by Hartogs lemma (1.6.2).

Let $(\f_j) \in {\mathcal F}_{\mu}^{\N}$. Then 
$\p_j:=\f_j-\sup_X \f_j \in {\mathcal F}_{0}$ which is relatively compact.
Assume first $\mu$ is smooth. 
Then $(\int_X \p_j d\mu)$ is bounded: this is because if
$\p_{j_k} \rightarrow \p$ in $L^1(X)$ then
$\p_{j_k} \mu \rightarrow \p \mu$ in the weak sense of (negative) measures
hence $\int_X \p_{j_k} d \mu \rightarrow \int_X \p d \mu>-\infty$.
Now $\int_X \p_j d\mu=\int_X \f_j d\mu -\int_X \sup_X \f_j d\mu=-\sup_X \f_j$ 
thus $(\sup_X \f_j)$ is bounded and we can apply the previous proposition
to conclude that $(\f_j)$ is relatively compact ( it cannot converge
uniformly to $-\infty$ since $\int_X \f_j d\mu=0$).
 
When $\mu$ is not smooth, it only remains to prove
that $(\int_X \p_j d\mu)$ is bounded. Assume on the contrary 
that $\int_X \p_j d\mu \rightarrow -\infty$. Extracting a subsequence
if necessary we can assume $\int_X \p_j d\mu \leq -2^j$. Set
$\p=\sum_{j \geq 1} 2^{-j} \p_j$. This is a decreasing sequence of 
$\o$-psh functions, hence $\p \in PSH(X,\o)$ or $ \p \equiv -\infty$.
Now it follows from the previous discussion that $\int_X \p_j dV \geq -C$
if $dV$ denotes some smooth probability measure on $X$. Thus
$\int_X \p dV>-\infty$ hence $\p \in PSH(X,\o)$. We obtain a contradiction
since by the Monotone convergence theorem,
$\int_X \p d\mu=\sum_{j \geq 1 } 2^{-j} \int_X \p_j d\mu=-\infty$.
\end{proof}

\begin{exa}
It was part of our definition 1.1 that $\o$-psh functions are integrable
with respect to a fixed volume form. Therefore
$PSH(X,\o) \subset L^1(\mu)$ for every smooth probability measure $\mu$ 
on $X$. More generally if $\mu$ is a probability measure on $X$ such
that
\begin{equation}
\mu=\Theta+dd^c(S),
\end{equation}
where $\Theta$ is smooth and $S$ is a {\bf positive} current of bidimension
$(1,1)$ on $X$, then $PSH(X,\o) \subset L^1(\mu)$ for any smooth
$\o$. Indeed let $\f$ in $PSH(X,\o)$, $\f \leq 0$. 
If $\f$ is smooth, it follows from Stokes theorem that
\begin{eqnarray*}
0 \leq \int_X (-\f) d\mu &=& \int_X (-\f) \Theta+\int_X (-\f) dd^c S \\
&\leq& C_{\Theta} ||\f||_{L^1}+\int_X S \wedge (-dd^c \f) \\
&\leq& C_{\Theta} ||\f||_{L^1}+\int_X S \wedge \o<+\infty,
\end{eqnarray*}
where the last inequality follows from $S \geq 0$ and 
$-dd^c\f \leq \o$. The general case follows by regularizing $\f$ (see
Appendix).

Probability measures satisfying $(1)$ naturally arise in complex dynamics
(see [25]).
Observe also that Monge-Amp\`ere measures arising from the local theory
of Bedford and Taylor [5] do satisfy $(1)$: if $u$ is psh and locally
bounded near e.g. the unit ball $B$ of $\C^n$, we can extend it to
$\C^n$ as a global psh function with logarithmic growth considering
$$
U(z):=\left\{ \begin{array}{rl}
u(z) & \text{ if } z \in B \\
\max(u(z),A \log^+|z|-\sup_B|u|-1) & 
\text{ if } z \in (1+\e)B \setminus B \\
A \log^+|z|-\sup_B|u|-1 & \text{ if } z \in \C^n \setminus (1+\e)B
\end{array} \right.
$$
where $\log^+|z|:=\max(\log|z|,0)$ and
with $A$ large enough. We assume $A=1$ for simplicity. Now
$\f:=U-\frac{1}{2}\log[1+|z|^2]+C$ extends as a bounded function in 
$PSH(\C\P^n,\o)$, where $\o$ is the Fubini-Study K\"ahler form on $\C\P^n$,
so $\f \geq 0$ if $C>0$ is large enough.
To conclude note that, setting $\o_{\f}:=\o+dd^c \f \geq 0$,
we get $\o_{\f}^n=(dd^c u)^n$ in $B$ and
$$
\o_{\f}^n=\o^n+dd^c S, \text{ where }
S=\f \sum_{j=0}^{n-1} \o_{\f}^j \wedge \o^{n-1-j} \geq 0.
$$
The Monge-Amp\`ere operator $\o_{\f}^n$ will be 
defined in the next section.
\end{exa}

\begin{exa}
If $\mu$ is a probability measure on $X=\C\P^n$ and $\o$ denotes as before the
Fubini-Study K\"ahler form, then
$$
\f_{\mu}(x):=\int_{\C\P^n} \log \left( 
\frac{||x \wedge y||}{||x|| \cdot ||y||} \right) d\mu(y)
$$
defines a $\o$-psh function on $\C\P^n$. Such functions have been considered
by Molzon, Shiffman and Sibony [34], [33] in order to define capacities
on $\C\P^n$. However they do not characterize pluripolar sets.
\end{exa}

\section{Monge-Amp\`ere capacity}

We assume in this section that $\o$ is a K\"ahler form on $X$. 
Let $T$ be a positive closed current of bidegree $(p,p)$ on $X$,
$0 \leq p \leq n=\dim_{\C} X$.
It can be thought of as a closed differential form of
bidegree $(p,p)$ with measure coefficients whose total variation is
controlled by
$$
||T||:=\int_X T \wedge \o^{n-p}.
$$ 
We refer the reader to chapter 3 of [17] 
for basic properties of  positive currents. 
Given $\f \in PSH(X,\o)$ we write $\f \in L^1(T)$ if $\f$ is integrable
with respect to each (measure) coefficient of $T$. This is equivalent
to $\f$ being integrable with respect to the trace measure 
$T \wedge \o^{n-p}$. In this case 
the current $\f T$ is well
defined, hence so is 
$$
\o_{\f} \wedge T:=\o \wedge T+dd^c( \f T).
$$
This is again a {\it positive} closed current on $X$, of bidegree $(p+1,p+1)$.
Indeed positivity is a local property which is stable under taking limits.
One can locally regularize $\f$ and approximate $\o_{\f} \wedge T$ by
the currents $\o_{\f_{\e}} \wedge T$ which are positive since
$\o_{\f_{\e}}$ are smooth positive forms. 

When $\f \in PSH(X,\o) \cap L^{\infty}(X)$ then $\f \in L^1(T)$ for any
positive closed current $T$ of bidegree $(p,p)$.
One can thus inductively define  
$\o_{\f}^j \wedge T$, $1 \leq j \leq n-p$, for 
$\f \in PSH(X,\o) \cap L^{\infty}(X)$. 
For $T=0$ and $j=n$ one obtains
the {\it complex Monge-Amp\`ere operator},
$\o_{\f}^n$.
It follows from the local theory that the operator $\f \mapsto \o_{\f}^n$
is continuous under monotone sequences (see [5]). The proof of these 
continuity properties is simpler in our compact setting. 
We refer the reader to [26] where this is proved 
in a more general global context.

\begin{pro}[Chern-Levine-Nirenberg inequalities]
Let $T$ be a positive closed current of bidegree $(p,p)$ on $X$
and $\f \in PSH(X,\o) \cap L^{\infty}(X)$.

Then
$||\o_{\f} \wedge T||=||T||$.
Moreover if $\p \in PSH(X,\o) \cap L^1(T)$, 
then $\p \in L^1(T \wedge \o_{\f})$ and 
$$
||\p||_{L^1(T \wedge \o_{\f})} \leq ||\p||_{L^1(T)}+
[2 \sup_X \p + \sup_X \f -\inf_X \f] ||T||.
$$
\end{pro}

\begin{proof}
By Stokes theorem, $\int_X dd^c \f \wedge T \wedge \o^{n-p-1}=0$, hence
$$
||\o_{\f} \wedge T||:=\int_X \o_{\f} \wedge T \wedge \o^{n-p-1}
=\int_X T \wedge \o^{n-p}:=||T||.
$$
Consider now $\p \in L^1(T)$. Since $T$ has measure coefficients, this
simply means that $\p$ is integrable with respect to the total variation
of these measures. Assume first $\p \leq 0$, $\f \geq 0$
and $\f,\p$ are smooth. Then
$$
||\p||_{L^1(T \wedge \o_{\f})}:= \int_X (-\p)T \wedge \o_{\f} \wedge \o^{n-p-1}
=||\p||_{L^1(T)}+
\int_X (-\p) T \wedge dd^c \f \wedge \o^{n-p-1}.
$$
Now it follows from Stokes theorem that
$$
\int_X (-\p) T \wedge dd^c \f \wedge \o^{n-p-1}=
\int_X \f T \wedge (-dd^c \p) \wedge \o^{n-p-1} 
$$
$$
\leq
\int_X \f T \wedge \o^{n-p} \leq \sup_X \f \int_X  T \wedge \o^{n-p},
$$
where the forelast inequality follows from $ \f T \wedge \o^{n-p} \geq 0$
and $-dd^c \p \leq \o$. This yields
$$
||\p||_{L^1(T \wedge \o_{\f})} \leq ||\p||_{L^1(T)}+\sup_X \f ||T||.
$$
The general case follows by regularizing $\f,\p$, 
observing that $\o_{\f}=\o_{\f'}$
where $\f'=\f-\inf_X \f \geq 0$, and decomposing
$\p=\p'+\sup_X \p$ with $\p'=\p-\sup_X \p \leq 0$. 
\end{proof}

\begin{rqe}
The fact that the $L^1$-norm of $\p$ with respect to the probability
measure $T \wedge \o_{\f} \wedge \o^{n-p-1}$ is controlled by its
$L^1$-norm with respect to $T \wedge \o^{n-p}$ is similar to the
phenomenon already encountered in example 1.8: one can write
$$
T \wedge \o_{\f} \wedge \o^{n-p-1}=T \wedge \o^{n-p}+dd^c S,
\; S=(\f-\inf_X \f) T \wedge \o^{n-p-1} \geq 0.
$$
This type of estimates is usually referred to as "Chern-Levine-Nirenberg
inequalities", in reference to [10] where simpler
-but fondamental- $L^{\infty}$-estimates were established 
(with $\p=constant$). 
Estimates involving the $L^1$-norm of $\p$ were first proved in the local
context by Cegrell [8] and Demailly [14].
\end{rqe}

A straightforward induction yields the following

\begin{cor}
Let $\p,\f \in PSH(X,\o)$ with $0 \leq\f \leq 1$. Then
$$
0 \leq \int_X |\p| \o_{\f}^n \leq \int_X |\p| \o^n+n[1+2\sup_X \p] \int_X
\o^n.
$$
\end{cor}

Following Bedford-Taylor [5] and Kolodziej [31] 
we introduce the following Monge-Amp\`ere capacity.

\begin{defi}
Let $K$ be a Borel subset of $X$. We set
$$
Cap_{\o}(K):=\sup \left\{ \int_K \o_{\f}^n \, / \, \f \in PSH(X,\o), 
0 \leq \f \leq 1 \right\}.
$$
\end{defi}

Note that this definition only makes sense when the cohomology 
class $[\o]$ is big, i.e. when $[\o]^n>0$, and when it admits
locally bounded potentials $\f$. This implies, by a regularization
result of Demailly [14], that $[\o]$ is big and nef. To simplify
we actually assume that $\o$ (hence $[\o]$) is K\"ahler.

\begin{pro} \text{ }

1) If $K \subset K' \subset X$ are Borel subsets then
$$
Vol_{\o}(K):=\int_K \o^n \leq Cap_{\o}(K) 
\leq Cap_{\o}(K') \leq Cap_{\o}(X)=Vol_{\o}(X).
$$

2) If $K_j$ are Borel subsets of $X$ then
$Cap_{\o}(\cup K_j) \leq \sum Cap_{\o}(K_j)$.
Moreover $Cap_{\o}(\cup K_j)=\lim Cap_{\o}(K_j)$ if $K_j \subset K_{j+1}$.

3) If $\o_1 \leq \o_2$ then $Cap_{\o_1}(\cdot) \leq Cap_{\o_2}(\cdot)$.
For all $A \geq 1$, 
$Cap_{\o}(\cdot) \leq Cap_{A \o}(\cdot) \leq A^n Cap_{\o}(\cdot)$.
In particular if $\o,\o'$ are two K\"ahler forms then there exists 
$C \geq 1$ such that
$$
\frac{1}{C}Cap_{\o}(\cdot) \leq Cap_{\o'}(\cdot) \leq C \cdot Cap_{\o}(\cdot).
$$

4) If $f:X \rightarrow X$ is holomorphic then for all Borel subset
$K$ of $X$,
$$
Cap_{\o}(f(K)) \leq Cap_{f^*\o}(K).
$$
In particular $Cap_{\o}(f(K)) = Cap_{\o}(K)$ for every 
$\o$-isometry $f$.
\end{pro}

\begin{proof}
That $Vol_{\o}(\cdot) \leq Cap_{\o}(\cdot)$ is a straightforward consequence
of the definition (since $\o_0 =\o$). It then follows from Stokes theorem
that $\int_X \o_{\f}^n=\int_X \o^n$ for every 
$\f \in PSH(X,\o) \cap L^{\infty}(X)$, thus 
$Vol_{\o}(X) = Cap_{\o}(X)$.

Property 2) is a straightforward consequence of the definitions.

If $\o_1 \leq \o_2$ then $PSH(X,\o_1) \subset PSH(X,\o_2)$ hence
$Cap_{\o_1}(\cdot) \leq Cap_{\o_2}(\cdot)$. Fix $A \geq 1$. If 
$\p \in PSH(X,A \o)$ is such that $0 \leq \p \leq 1$ then
$\p/A \in PSH(X,\o)$ with $0 \leq \p/A \leq 1/A \leq 1$.
Moreover $(A \o +dd^c \p)^n=A^n(\o+dd^c(\p/A))^n$. This shows
$Cap_{A \o}(\cdot) \leq A^n Cap_{\o}(\cdot)$. 

In particular if $\o,\o'$ are both K\"ahler then 
$A^{-1} \o \leq \o' \leq A \o$ for some constant $A \geq 1$, hence 
$C^{-1} Cap_{\o}(\cdot) \leq Cap_{\o'}(\cdot) \leq C \cdot Cap_{\o}(\cdot)$
with $C=A^n$.

It remains to prove 4). It follows from the change of variables formula that
if $\f \in PSH(X,\o)$ with $0 \leq \f \leq 1$ then
$$
\int_{f(K)} \o_{\f}^n \leq \int_K f^* \o_{\f}^n=
\int_K (f^*\o+dd^c(\f \circ f))^n
\leq Cap_{f^*\o}(K)
$$
since $\f \circ f \in PSH(X,f^* \o)$ with $0 \leq \f \circ f \leq 1$.
We infer $Cap_{\o}(f(K)) \leq Cap_{f^*\o}(K)$.
When $f$ is a $\o$-isometry, i.e. $f \in Aut(X)$ with $f^*\o=\o$, then
the mapping $\f \mapsto \f \circ f$ is an isomorphism of
$\{u \in PSH(X,\o) \, / \, 0 \leq u \leq 1 \}$, whence
$Cap_{\o}(f(K))=Cap_{f^*\o}(K)=Cap_{\o}(K)$.
\end{proof}

Let $PSH^-(X,\o)$ denote the set of negative $\o$-psh functions.
A set is said to be $PSH(X,\o)$-{\it polar} if it is included in the 
$-\infty$ locus of some function $\p \in PSH(X,\o)$, $\p \not\equiv -\infty$.
As we shall soon see, the sets of zero Monge-Amp\`ere capacity are
precisely the $PSH(X,\o)$-polar sets. We start by establishing the following:

\begin{pro}
If $P$ is a $PSH(X,\o)$-polar set, then $Cap_{\o}(P)=0$.
More precisely, if $\p \in PSH^-(X,\o)$ then
$$
Cap_{\o}(\p<-t) \leq \frac{1}{t}\left[\int_X (-\p) \o^n+n Vol_{\o}(X) \right],
\; \forall t>0.
$$
\end{pro}

\begin{proof}
Fix $\f \in PSH(X,\o)$ such that $0 \leq \f \leq 1$. Fix $t>0$ and set
$K_t=\{x \in X \, / \, \p(x)<-t \}$. By Chebyshev's inequality,
$$
\int_{K_t} \o_{\f}^n \leq \int_X (-\p/t) \o_{\f}^n \leq
\frac{1}{t}\left[\int_X (-\p) \o^n+n Vol_{\o}(X) \right],
$$
where the last inequality follows from corollary 2.3.
Taking supremum over all $\f'$s yields the claim.
\end{proof}

Observe that the previous proposition says that $Cap_{\o}^*(P)=0$,
where
$$
Cap_{\o}^*(E):=\inf \{Cap_{\o}(G) \, / G \text{ open with } E \subset G \},
$$
is the outer capacity associate to $Cap_{\o}$.

Our aim is now to show that $\o$-psh functions are quasicontinuous
with respect to $Cap_{\o}$ (corollary 2.8). We first need to show
that decreasing sequences of $\o$-psh functions converge
"in capacity".

\begin{pro}
Let $\p,\p_j \in PSH(X,\o) \cap L^{\infty}(X)$ such that
$(\p_j)$ decreases to $\p$. 
Then for each $\d>0$,
$$
Cap_{\o}(\{\p_j>\p+\d\}) \rightarrow 0.
$$
\end{pro}

\begin{proof}
We can assume w.l.o.g. that $Vol_{\o}(X)=1$ and $0 \leq \p_j-\p \leq 1$.
Fix $\d>0$ and $\f \in PSH(X,\o)$, $0 \leq \f \leq 1$.
By Chebyshev inequality, it suffices to control 
$\int_X (\p_j-\p) \o_{\f}^n$ uniformly in $\f$. It follows from Stokes 
theorem that
$$
\int_X (\p_j-\p) \o_{\f}^n=\int_X (\p_j-\p) \o \wedge \o_{\f}^{n-1}
-\int_X d(\p_j-\p) \wedge d^c \f \wedge \o_{\f}^{n-1}.
$$
Now by Cauchy-Schwartz inequality,
$$
\left|\int_X d f_j \wedge d^c \f \wedge \o_{\f}^{n-1} \right|
\leq \left( \int_X d f_j \wedge d^c f_j
\wedge \o_{\f}^{n-1} \right)^{1/2} \cdot
\left( \int_X d \f \wedge d^c \f 
\wedge \o_{\f}^{n-1} \right)^{1/2},
$$
where we set $f_j:=\p_j-\p \geq 0$.
Moreover
$$
\int_X d \f \wedge d^c \f \wedge \o_{\f}^{n-1}=
\int_X \f (-dd^c \f) \wedge \o_{\f}^{n-1} \leq 
\int_X \f \o \wedge \o_{\f}^{n-1} \leq 1,
$$
since $\f \o_{\f}^{n-1} \geq 0$, $-dd^c \f \leq \o$ and $\f \leq 1$.
Similarly
$$
\int_X df_j \wedge d^c f_j \wedge \o_{\f}^{n-1}=
\int_X -f_jdd^c f_j \wedge \o_{\f}^{n-1}
\leq \int_X f_j \o_{\p} \wedge \o_{\f}^{n-1}.
$$
Altogether this yields
\begin{eqnarray*}
\int_X (\p_j-\p) \o_{\f}^n &\leq& \int_X (\p_j-\p) \o \wedge \o_{\f}^{n-1}
+\left( \int_X (\p_j-\p) \o_{\p} \wedge \o_{\f}^{n-1} \right)^{1/2} \\
&\leq& \sqrt{2} \left( \int_X (\p_j-\p) (\o+\o_{\p}) 
\wedge \o_{\f}^{n-1} \right)^{1/2},
\end{eqnarray*}
where the last inequality follows from the elementary inequalities
$0 \leq a \leq \sqrt{a} \leq 1$ and 
$\sqrt{a}+\sqrt{b} \leq \sqrt{2} \sqrt{a+b}$.

Going on replacing at each step a term $\o_{\f}$ by $\o+\o_{\p}$,
we end up with
$$
\int_X (\p_j-\p) \o_{\f}^n \leq 2 \left( 
\int_X (\p_j-\p) (\o+\o_{\p})^n \right)^{1/2^n}.
$$
The majorant being independent of $\f$ and converging to $0$ as
$j \rightarrow +\infty$ (by dominated convergence theorem), this 
completes the proof.
\end{proof}

\begin{cor}[Quasicontinuity]
Let $\f \in PSH(X,\o)$. For each $\e>0$ there exists an open subset
$O_{\e}$ of $X$ such that $Cap_{\o}(O_{\e})<\e$
and $\f$ is continuous on $X \setminus O_{\e}$.
\end{cor}

\begin{proof}
For $t>0$ large enough, the set $O_1=\{\f<-t\}$ has capacity $<\e/2$
by proposition 2.6. Working in $X \setminus O_1$ we can thus replace
$\f$ by $\f_t=\max(\f,-t)$ which is bounded on $X$.
Regularizing $\f$ (see Appendix), we can find a sequence $\p_j$ of
smooth $A\o$-psh functions which decrease to $\f_t$ on $X$,
for some $A \geq 1$.
By proposition 2.7, the set $O_j=\{ \p_{k_j} >\f_t+1/j \}$ has capacity
$<\e2^{-j-1}$ if $k_j$ is large enough.
Now $\p_{k_j}$ uniformly converges to $\f=\f_t$ on
$X \setminus O_{\e}$, $O_{\e}=\cup_{j \geq 1} O_j$, so
$\f$ is continuous on $X \setminus O_{\e}$ and 
$Cap_{\o}(O_{\e}) \leq \e$.
\end{proof}

\begin{exa}
The capacity $Cap_{\o}(\cdot)$ does not distinguish between
"big sets". Assume indeed there exists an ample divisor $D$
such that $[D] \sim k \o$, $k \in \N$. Then there exists
$\f \in PSH(X,\o)$ such that $dd^c\f=k^{-1}[D]-\o$. Note that
$\f \in {\mathcal C}^{\infty}(X \setminus D)$, 
$e^{\f} \in {\mathcal C}^0(X)$ and $\{\f=-\infty\}=D$.
Replacing $\f$ by $\f-\sup_X \f$ if necessary, we may assume
$\sup_X \f=0$.
Consider $\f_c=\max(\f,-c) \in PSH(X,\o) \cap {\mathcal C}^0(X)$.
Then $\f_c \equiv \f$ outside some neighborhood $V_c=\{\f <-c \}$ of $D$.
Since $0 \leq 1+\f_1 \leq 1$ and $\o_{1+\f_1}=\o_{\f}=0$ in 
$X \setminus V_1$, we get
$$
Cap_{\o}(X)=\int_X (\o_{1+\f_1})^n=\int_{V_1} (\o_{1+\f_1})^n
\leq Cap_{\o}(V_1),
$$
hence $Cap_{\o}(V_1)=Cap_{\o}(X)$.

As a concrete example take $X=\C\P^n$ and  $\o=\o_{FS}$, 
$D$ being some hyperplane $H_{\infty}$ "at infinity" ($k=1$).
Set $\f[z:t]=\log|t|-\frac{1}{2}\log[||z||^2+|t|^2]$ where
$z$ denotes the euclidean coordinates in $\C^n=\C\P^n \setminus H_{\infty}$
and $H_{\infty}=(t=0)$. Observe that $\sup_{\C\P^n} \f=0$.
One then computes
$$
\C\P^n \setminus V_1=\left\{ z \in \C^n \, / \, |z| \leq \sqrt{e^2-1} \right\}.
$$
Thus the capacity of the complement of
any euclidean ball of radius smaller than
$\sqrt{e^2-1}$ equals $1$.
\end{exa}

The definition of $Cap_{\o}$ mimics the definition of the relative
Monge-Amp\`ere capacity introduced by Bedford and Taylor in [5].
Fix ${\mathcal U}=\{ {\mathcal U}_{\a} \}$ a finite covering of
$X$ by strictly pseudoconvex open subsets of $X$, 
${\mathcal U}_{\a}=\{x \in X \, / \, \r_{\a}(x) <0 \}$, where
$\r_{\a}$ is a strictly psh smooth function defined in a neighborhood
of $\overline{{\mathcal U}_{\a}}$. Fix $\d>0$ such that
${\mathcal U}^{\d}=\{ {\mathcal U}_{\a}^{\d} \}$ is still a covering
of $X$, where ${\mathcal U}_{\a}^{\d}=\{x \in X \, / \, \r_{\a}(x)<-\d \}$.
For a Borel subset $K$ of $X$, we set
$$
Cap_{BT}(K):=\sum_{\a} Cap_{BT}(K \cap {\mathcal U}_{\a}^{\d}, {\mathcal
  U}_{\a}),
$$
where
$$
Cap_{BT}(E, \Omega):=\sup \left\{ \int_E (dd^c u)^n \, / \,
u \in PSH(\Omega), \, 0 \leq u \leq 1 \right\}
$$
is the capacity studied by Bedford and Taylor.
The next proposition is due to Kolodziej [31]. We include a slightly
different proof.

\begin{pro}
There exists $C \geq 1$ such that
$$
\frac{1}{C} Cap_{\o}(\cdot) \leq Cap_{BT}(\cdot) \leq C \cdot
Cap_{\o}(\cdot).
$$
\end{pro}

\begin{proof}
Let $E$ be a Borel subset of $X$. Since
$Cap_{\o}(E \cap {\mathcal U}_{\a}^{\d}) \leq Cap_{\o}(E)
\leq \sum_{\a} Cap_{\o}(E \cap {\mathcal U}_{\a}^{\d})$, it is sufficient
to show that if $\Omega=\{x \in X \, / \, \r(x)<0 \}$ is a smooth
hyperconvex subset of $X$, then there exists $C \geq 1$ such that
for all $E \subset \Omega_{\d}$,
$$
\frac{1}{C} Cap_{\o}(E)  \leq Cap_{BT}(E,\Omega)  \leq C \cdot 
Cap_{\o}(E),
$$
where $\Omega_{\d}=\{x \in X \, / \, \r(x)<-\d \}$.

It is an easy and well known fact in the local theory that the capacities
$Cap(\cdot,\Omega)$ and $Cap(\cdot,\Omega')$ are comparable when 
$\Omega' \subset \Omega$ (see e.g. theorem 6.5 in [14]).
Therefore we can assume (passing to a finer covering if necessary) that
$\o=dd^c \p$ near $\overline{\Omega}$. Fix $C_1>0$ such that
$-C_1\leq \p \leq C_1$ on $\Omega$. Fix $\f \in PSH(X,\o)$ such that
$0 \leq \f \leq 1$ on $X$ and set
$u=(2C_1)^{-1}(\f+\p+C_1)$. Then $u \in PSH(\Omega)$ and 
$0 \leq u \leq 1$, hence
$$
\int_E \o_{\f}^n=(2C_1)^n \int_E (dd^c u)^n \leq (2 C_1)^n Cap_{BT}(E,\Omega),
$$
which yields $Cap_{\o}(E) \leq (2C_1)^n  Cap_{BT}(E,\Omega)$.
Observe that we have not used here that $\o$ is K\"ahler.

For the reverse inequality we consider $\chi \in {\mathcal C}^{\infty}(X)$
such that $\chi \equiv 0$ in $X \setminus \Omega$ and 
$\chi<0$ in $\Omega$.
Replacing $\chi$ by $\e \chi$ if necessary, we can assume $\chi \in
PSH(X,\o)$.
This is because $\o$ is K\"ahler (and this is the only place where we shall 
use this crucial assumption). Fix $\e>0$ so small that $\chi \leq -\e$
on $\Omega_{\d}$. Let now $u \in PSH(\Omega)$ be such that
$0 \leq u \leq 1$ on $\Omega$. Consider
$$
\f(x)=\left\{ \begin{array}{cl}
\frac{u-\p+C_1}{2+2C_1} & \text{ in } \Omega_{\d} \\
\max\left(\frac{u-\p+C_1}{2+2C_1}, \frac{2}{\e}\chi(x)+1\right) 
& \text{ in } \Omega \setminus \Omega_{\d} \\
1 & \text{ in } X \setminus \Omega
\end{array} \right.
$$
Observe that $0 \leq u':=(u-\p+C_1)/(2+2C_1) \leq (1+2C_1)/(2+2C_1)<1$
in $\Omega$. Therefore $\f \in PSH(X,\frac{2}{\e}\o)$ since
$\frac{2}{\e}\chi(x)+1 \leq -1 <u'$ in $\Omega_{\d}$,
while
$\frac{2}{\e}\chi(x)+1 \equiv 1>u'$ on $\partial \Omega$.
Note also that $0 \leq \f \leq 1$ thus for $E \subset \Omega_{\d}$,
\begin{eqnarray*}
\frac{1}{(2+2C_1)^n}\int_E (dd^c u)^n &=&
\int_E \left( \frac{\o}{2+2C_1}+dd^c \f \right)^n \leq
\int_E \left( \frac{2}{\e}\o+dd^c \f \right)^n \\
&\leq&
Cap_{2\o/\e}(E) \leq \left(\frac{2}{\e}\right)^n Cap_{\o}(E)
\end{eqnarray*}
hence $Cap_{BT}(E,\Omega)  \leq 4^n(1+C_1)^n\e^{-n} Cap_{\o}(E)$.
\end{proof}

Since locally pluripolar sets are precisely
the sets of zero relative capacity [5], we obtain the following

\begin{cor}
$Cap^*_{\o}(P)=0 \Leftrightarrow Cap^*_{BT}(P)=0 \Leftrightarrow P 
\text{ is locally pluripolar}.$
\end{cor}

We shall show later on that  locally pluripolar sets
are $PSH(X,\o)$-polar when $\o$ is K\"ahler (see theorem 6.2).

The following two results are 
direct consequences of the corresponding
results of Bedford and Taylor [4], [5].

\begin{thm}[Dirichlet Problem]
Let $\f \in PSH(X,\o) \cap L^{\infty}(X)$. Let $B$ be a small ball
in $X$. Then there exists $\hat{\f} \in PSH(X,\o)$ such that
$\hat{\f}=\f$ in $X \setminus B$, $\hat{\f} \geq \f$
and $(\o_{\hat{\f}})^n=0$ in $B$.
Moreover if $\f_1 \leq \f_2$ then $\hat{\f_1} \leq \hat{\f_2}$.
\end{thm}

\begin{thm}[Comparison principle]
Let $\f,\p \in PSH(X,\o) \cap L^{\infty}(X)$. Then
$$
\int_{\{\f<\p\}} \o_{\p}^n \leq \int_{\{\f <\p\}} \o_{\f}^n.
$$
\end{thm}

\section{The relative extremal function}

We now introduce a substitute for the relative extremal function which has
revealed so useful in the local theory [5]. Namely if $E$ is a Borel subset of
$X$, we set
$$
h_{E,\o}(x):=\sup \left\{\f(x) \, / \, \f \in PSH(X,\o) , \, 
\f \leq 0 \text{ and } \f_{|E} \leq -1 \right\}.
$$
We let $h_{E,\o}^*$ denote its upper-semi-continuous regularization, which we 
call the relative $\o$-plurisubharmonic extremal function of the subset 
$E \subset X$. It enjoys several natural properties; we list some of them
below. The proofs follow from standard arguments together with theorem 6.2:
\begin{itemize}
\item The function $h_{E,\o}^*$ is $\o$-psh. It satisfies
$-1 \leq h_{E,\o}^* \leq 0$ on $X$
and $h_{E,\o}^*=-1$ on $E \setminus P$, where $P$ is pluripolar.

\item If $E \subset X$ and $P \subset X$ is pluripolar, then 
$h_{E \setminus P}^* \equiv h_E^*$; in particular $h_{P}^* \equiv 0$.

\item If $(E_j)$ increases towards $E \subset X$, then
$h_{E_j}^*$ decreases towards $h_E^*$.

\item If $(K_j)$ is a sequence of compact subsets decreasing towards
$K$, then $h_{K_j}^*$ increases (a.e.) towards $h_K^*$.
\end{itemize}

As in the local theory, the complex Monge-Amp\`ere of the relative extremal
function of a subset $E \subset X$ vanishes outside $\overline{E}$, except
perhaps on the set $\{h_E^*=0\}$ which, in the local theory, lies in
the boundary of the domain.

\begin{pro} \text{ }
When the open set $\Omega_E:=\{x \in X \, / \, h_{E,\o}^*(x)<0 \}$ is
   non-empty, then
$$
\left( \o_{h_{E,\o}^*} \right)^n=\left(\o+dd^c h_{E,\o}^* \right)^n
=0 \text{ in } \Omega_E \setminus \overline{E}.
$$
\end{pro}

\begin{proof}
Assume  that $\Omega_E:=\{h_{E,\o}^*<0\}$ is non-empty.
It follows from Choquet's lemma that there exists an increasing sequence
$\f_j$ of $\o$-psh functions such that $\f_j=-1$ on $E$, $\f_j \leq 0$ on $X$,
and $h_{E,\o}^*=(\lim \f_j)^*$.
Let $a \in \Omega_E \setminus \overline{E}$ 
and fix a small ball $B \subset \Omega_E$ centered at 
point $a$. Let $\widehat{\f_j} = \widehat{(\f_j)}_B$ denote the functions obtained by applying
theorem 2.12, so that $(\o_{\widehat{\f_j}})^n=0$ in $B$.
If $B$ is chosen small enough, then $\widehat{\f_j} <0$ in $B$, hence
$\widehat{\f_j} \leq 0$ on $X$, while $\widehat{\f_j}=\f_j =-1$ on $E$.
This can be seen by showing that $\widehat{0}_B \rightarrow 0$ as the radius of
the ball $B$ shrinks to $0$.
Therefore $\lim \nearrow \widehat{\f_j}=h_{E,\o}^*$, hence 
$(\o_{h_{E,\o}^*})^n \equiv 0$
in a neighborhood of $a$, which prove our claim.
\end{proof}

We now establish an important result which expresses the capacity
in terms of the relative extremal function for any subset. It will
show in particular
that the set function ${Cap}^*_{\o}$ is a capacity in the sense of
Choquet which is outer regular.
For simplicity, we write $h_E$ for $h_{E,\o}$.

\begin{thm} 
Let $E \subset X$ be any Borel subset, then
$$
{Cap}_{\o}^* (E) =  \int_X (-  h_{E}^*) \o_{h_{E}^*}^n.
\hskip2cm (\dag)
$$
The Monge-Amp\`ere capacity satisfies the following continuity properties:

 1) If $(E_{j})_{j \geq 0}$ is an increasing sequence 
 of arbitrary subsets of $X$ and $E := \cup_{j \geq 0} E_j$ then
$$
{Cap}_{\o}^* (E) =  \lim_{j \to + \infty}
 {Cap}_{\o}^* (E_{j}).
$$

 2) If $(K_{j})_{j \geq 0}$ is a decreasing sequence 
 of compact subsets of $X$ and $K := \cap_{j \geq 0} K_j$ then
 $$
{Cap}_{\o} (K) = {Cap}_{\o}^* (K) =  \lim_{j \to + \infty}
 {Cap}_{\o} (K_{j}).
$$

In particular ${Cap}_{\o}^*(\cdot)$ is an outer regular Choquet 
capacity on $X$.
\end{thm}

\begin{proof} 
We first establish $(\dag)$ when
$E = K \subset X$ is compact.
 Observe that in the definition of the capacity, it 
is enough to restrict ourselves to $\o$-psh functions 
$\f$ such that $- 1 < \f < 0$. Let $\f$ be such a function.  
The pluripolar set $N := \{h_K < h_K^*\}$ is of measure $0$ for 
the measure $\o_{\f}^n$. Since $- 1 < \f$, we  
have $K \subset N \cup \{h_K^* < \f \}$, hence 
the comparison principle yields 
$$ 
\int_{K} \o_{\f}^n \leq  \int_{\{h_K^* < \f \}} \o_{\f}^n \leq
 \int_{\{h_K^* < \f\}} \o_{h_E^*}^n .
$$
This shows
 ${Cap}_{\o} (K) \leq \int_{\Omega_K} \o_{h_K^*}^n$ since 
 $\{h_K^* < \f\} \subset \Omega_K$. 
It follows from proposition 3.1 that
$$
{Cap}_\o (K) \leq \int_{\Omega_K} (\o_{h_{K}}^*)^n = 
\int_{K}  (\o_{h_{K}^*})^n,
$$
whence equality.
 
Recall that $K \cap \{h_{K}^{*} > - 1 \} \subset \{h_{K} < h_{K}^{*}\} $ 
is of measure $0$ for $\o_{h_{K}^*}^n$, so 
$$
{Cap}_{\o}(K) = \int_{K} (\o_{h_{K}}^*)^n =
 \int_{K} (- h_{K}^*) (\o_{h_{K}}^*)^n =
 \int_{X} (- h_{K}^*) (\o_{h_{K}}^*)^n,
$$
 where the last equality follows from the fact that the equilibrium 
measure $(\o_{h_{K}^*})^n$ is supported on $K \cup \{h_K^* = 0\}.$
\vskip.1cm
 
  We assume now 
that $E = G \subset X$ is an {\it open} subset. Let $(K_j)$
be an exhaustive 
sequence of compact subsets of $G$ which increases to $G$. 
Since $h_{K_j}^* \downarrow h_G$ on $X$,
it follows from classical convergence results (see [40]) that
$(- h_{K_j}^*) \o_{h_{K_j}^*}^n \to (- h_{G}^*) \o_{ h_{G}}^n$ 
in the weak sense of measures on $X$, therefore
$$ 
\int_{X} (- h_{G}) \o_{h_{G}}^n =
  \lim_{j \to + \infty} \int_X (- h_{K_j}^*)  \o_{h_{K_j}^*}^n
=\lim_{j \rightarrow +\infty} Cap_{\o}(K_j)=Cap_{\o}(G),
$$
This proves $(\dag)$ when 
$E \subset X$ is an open subset.
\vskip.1cm

Finally let $E \subset X$ be any subset. By definition of the outer 
capacity, there is a sequence of open subsets $(O_j)_{j \geq 1}$ of $X$ containing $E$ such that
${Cap}_{\o}^* (E) = \lim_{j \to + \infty} {Cap}_{\o} (O_j).$
We can assume w.l.o.g. that the sequence $(O_j)_{j \geq 1}$ is decreasing.
  
  By a classical topological lemma of Choquet, there exists an 
increasing sequence $(u_j)_{j \geq 1}$  negative  $\o$-psh
functions on $X$ s.t. 
 $ u_j = - 1$ on $E$ with $u_j \uparrow h_E^*$ almost everywhere on $X$. 
We set for each $j \in \N,$ 
$G_j :=  O_j \cap \{ u_j < - 1 + 1 \slash j\}$. 
Then $(G_j)$ is a decreasing sequence of open subsets of $X$ such that $E
\subset G_j \subset O_j$ and $u_j - 1 \slash j \leq h_{G_j} \leq h_E$, 
so $h_{G_j} \uparrow h_E^*$ almost everywhere on $X.$ 
We infer $(- h_{G_j}) \o_{h_{G_j}}^n \to (- h_{E}^*) \o_{h_{E}}^n $ 
in the weak sense of measures on $X$, thus 
  $$ \int_{X} (- h_{E}^*) \o_{h_{E}^*}^n  =
  \lim_{j \to + \infty} \int_X  (- h_{G_j}^*) \o_{h_{G_j}}^n.$$
  On the other hand we have by construction 
${Cap}_{\o}^* (E) \leq \lim_{j \to + \infty} {Cap}_{\o}^* (G_j) \leq 
  \lim_{j \to + \infty} {Cap}_{\o}^* (O_j) = {Cap}_{\o}^* (E).$
  Therefore using $(\dag)$ for open subsets,  we get
  $$\int_{X} (- h_{E}^*) \o_{h_{E}^*}^n =  {Cap}_{\o}^* (E).$$
\vskip.1cm

Observe that 1) follows straightforwardly from this formula. Indeed if
$(E_j)$ increases towards $E$, then 
$h_{E_j}^*$ decreases towards $h_E^*$, hence
$(-h_{E_j}^*) (\o_{h_{E_j}^*})^n \rightarrow (-h_{E}^*) (\o_{h_{E}^*})^n$,
so that
$$
{Cap}_{\o}^* (E)=
\int_{X} (- h_{E}^*) \o_{h_{E}^*}^n=
\lim \int_{X} (-h_{E_j}^*) (\o_{h_{E_j}^*})^n=
\lim {Cap}_{\o}^* (E_j).
$$

It remains to prove 2). Let $(K_{j})$ be a decreasing 
sequence of compact subsets of $X$ which converges to $K.$ 
We claim that $h_{K_{j}}^{*} \uparrow h_K^*$ almost everywhere on $X$. Indeed, 
the extremal function $h_{K_{j}}^{*}$ 
increases almost everywhere 
to a $\o-$psh function $h$ such that $h \leq h_K^*$ on $X$.
We want to prove that $h_K \leq h$ on $X$. Let $u \in PSH (X,\o)$ such
that $u \leq 0$ on $X$ and $u_{|K} \leq -1$. 
Fix $\e > 0$ and consider the open
subset $G_{\e} := \{ u < - 1 + \e\}$. 
Then $K \subset G_{\e}$, thus 
$K_j \subset G_{\e}$ for $j$ large enough.
This yields $u - \e \leq h_{K_j}$ for $j$ large enough,
hence $u \leq h$ on $X$. Therefore $h_K \leq h$ on $X$ as claimed. 

Since
$(-h_{K_j})  \o_{h_{K_j}^*}^n $
converges weakly to $(-h_K^*)  \o_{h_{K}^{*}}^n$ on 
$X$,
we infer
$$
Cap_{\o}(K)=
\int_{X} (- h_{K}^*)  \o_{h_{K}^*}^n  =
\lim_{j \to + \infty} \int_{X} 
(- h_{K_j}^*)  \o_{h_{K_j}^*}^n
=\lim_{j \rightarrow +\infty} Cap_{\o}(K_j).
$$

From this last property, taking a decreasing sequence $(K_j)_{j \geq 0}$ of
compact subsets such that $K = \cap_{j \geq 0} K_j$ and $K_{j + 1} \subset
K_j^{\circ}$ for any $j \in \N$ we obtain
${Cap}_{\o}^* (K) = \lim_{j \to + \infty} {Cap}_{\o} (K_j) = {Cap}_{\o} (K).$ 
\end{proof}
\begin{cor} For any subset $P \subset X$, we have
$ Cap^*_{\o}(P)=0 \Leftrightarrow h_P^* \equiv 0. 
$
\end{cor}

\section{Alexander capacity}

We now introduce another capacity which is defined by means of a
global extremal function. It is closely related to the projective
capacity introduced by Alexander in [1].
We assume throughout this section that $\o$ is a closed real current on $X$
with continuous local potentials.

\subsection{Global extremal functions}

\begin{defi}
Let $K$ be a Borel subset of $X$. We set
$$
V_{K,\o}:=\sup \left\{ \f(x) \, / \, \f \in PSH(X,\o), \, \f \leq 0 
\text{ on } K \right\}.
$$
\end{defi}

This definition mimics the definition of the so-called "Siciak's extremal
function" usually defined for Borel subset of $X=\C\P^n$
that are bounded in $\C^n=X \setminus H_{\infty}$, where
$H_{\infty}$ denotes some hyperplane at infinity.
This function was introduced and studied by Siciak in [37],[38]
(see also [41]).
One can indeed check that this definition coincides with the classical
one if one chooses $\o=[H_{\infty}]$ to be the current of integration
along the hyperplane $H_{\infty}$. 
Similarly one could consider the case where $\o=[D]$ is the current of
integration along a positive divisor $D$ on $X$ and let $D$ play the role of
infinity. This approach has been used by some authors working in Arakelov
geometry to define capacities on projective varieties (see [32],[11]
and references therein). However this forces them to consider only compact
subsets of $X \setminus D$ and leads to less intrinsic notions of capacities.

In this article we always assume that the currents $\o$ involved
admit continuous potentials. 
This insures that the Monge-Amp\`ere operator $\o_{\f}^n$ is well-defined
on extremal functions $V_{K,\o}$.
If $\f \in L^1(X)$, we shall denote
by $\f^*$ its upper-semi-continuous regularization.

\begin{thm}
Let $K$ be a Borel subset of $X$.

1) $K$ is $PSH(X,\o)$-polar iff $\sup_X V_{K,\o}^* =+\infty$
iff $V_{K,\o}^* \equiv +\infty$.

2) If $K$ is not $PSH(X,\o)$-polar then $V_{K,\o}^* \in PSH(X,\o)$ and 
satisfies $V_{K,\o}^* \equiv 0$ in the interior of ${K}$,
$(\o_{V_{K,\o}^*})^n=0$ in $X \setminus \overline{K}$ and
$$
\int_{\overline{K}} (\o_{V_{K,\o}^*})^n=\int_X \o^n=Vol_{\o}(X).
$$
\end{thm}

\begin{proof}
Assume $\sup_X V_{K,\o}^*=+\infty$. By a lemma of Choquet
(see lemma 4.23 in [17], chapter 1),
we can find an increasing sequence of functions $\f_j \in PSH(X,\o)$
such that $\f_j =0$ on $K$ and $V_{K,\o}^*=(\lim\nearrow \f_j)^*$.
Extracting a subsequence if necessary, we can assume $\sup_X \f_j \geq 2^j$.
Set $\p_j=\f_j -\sup_X \f_j$. These functions belong
to  ${\mathcal F}_0$ which is a compact subfamily of $PSH(X,\o)$
(corollary 1.7). 
Recall that if $\mu$ is a smooth volume form on $X$ 
then there exists $C_{\mu}$ such that $\int \p_j d\mu \geq -C_{\mu}$
for all $j$.
Set $\p:=\sum_{j \geq 1} 2^{-j} \p_j$. Then 
$\p \in PSH(X,\o)$ as  a decreasing limit of functions in $PSH(X,\o)$
with $\int_X \p d\mu \geq -C_{\mu}>-\infty$.
Now for every $x \in K$ we get 
$\p(x)=-\sum_{j \geq 1} 2^{-j} \sup_X \f_j=-\infty$ hence 
$K \subset \{\p=-\infty\}$, i.e. $K$ is $PSH(X,\o)$-polar.

Conversely assume $K$ is $PSH(X,\o)$-polar, $K \subset \{\p=-\infty\}$
for some $\p \in PSH(X,\o)$. Then for all $c \in \R$, $\p+c \in PSH(X,\o)$
and $\p+c \leq 0$ on $K$. Therefore $V_{K,\o} \geq \p+c$, 
$\forall c \in \R$. This yields $V_{K,\o}=+\infty$ on
$X \setminus \{\p=-\infty\}$ hence $V_{K,\o}^* \equiv +\infty$ on
$X$ since $\{\p=-\infty\}$ has zero volume.
We have thus shown the following circle of implications:
$K \text{ is } PSH(X,\o)-\text{polar} \Rightarrow V_{K,\o}^* \equiv +\infty
\Rightarrow\sup_X V_{K,\o}^*=+\infty \Rightarrow K \text{ is }
PSH(X,\o)-\text{polar}$.  

Assume now that $K$ is not $PSH(X,\o)$-polar. Then 
$V_{K,\o}^* \in PSH(X,\o)$ (see proposition 1.6.2) and clearly satisfies
$V_{K,\o}^*=0$ in the interior of ${K}$. 
If we show that $(\o_{V_{K,\o}^*})^n=0$
in $X \setminus \overline{K}$ then
$$
\int_{\overline{K}} (\o_{V_{K,\o}^*})^n=\int_X (\o_{V_{K,\o}^*})^n=
\int_X \o^n,
$$
as follows from Stokes theorem. Let $\f_j \in PSH(X,\o)$ be an increasing 
sequence such that $\f_j=0$ on $K$ and $V_{K,\o}^*=(\lim\nearrow \f_j)^*$.
Fix $B$ a small ball in $X \setminus \overline{K}$.
Let $\hat{\f_j}$ be the solution of the Dirichlet problem with
boundary values $\f_j$. Then $\hat{\f_j} \in PSH(X,\o)$, 
$\hat{\f_j}=\f_j$ in $X \setminus B$ (in particular $\hat{\f_j}=0$
on $K$ hence $\hat{\f_j} \leq V_{K,\o}$) and the sequence
$(\hat{\f_j})$ is again increasing (theorem 2.12).
Since $(\o_{\hat{\f_j}})^n=0$ in $B$ and 
$(\lim\nearrow \hat{\f_j})=V_{K,\o}^*$, it follows from the continuity
of the complex Monge-Amp\`ere on increasing sequences
that $(\o_{V_{K,\o}^*})^n=0$ in $B$. As $B$ was an arbitrarily small
ball in $X \setminus \overline{K}$ we infer $(\o_{V_{K,\o}^*})^n=0$
in $X \setminus \overline{K}$.
\end{proof}

The following corollary has to be related to proposition 1.7.

\begin{cor}
Let $K$ be a Borel subset of $X$ and set
$$
{\mathcal F}_K:=\{\f \in PSH(X,\o) \, / \, \sup_K \f=0 \}.
$$
Then ${\mathcal F}_K$ is relatively compact iff $K$ is not $PSH(X,\o)$-polar.
\end{cor}

\begin{proof}
Observe that $V_{K,\o}(x)=\sup \{ \f(x) \, / \, \f \in {\mathcal F}_K \}$.
Thus if ${\mathcal F}_K$ is relatively compact then it is
uniformly bounded from above, hence $\sup_X V_{K,\o}<+\infty$, i.e.
$K$ is not $PSH(X,\o)$-polar.

Assume conversely that $K$ is not $PSH(X,\o)$-polar.
Let $(\f_j) \in {\mathcal F}_K^{\N}$. Then 
$\f_j \leq V_{K,\o} \leq \sup_X V_{K,\o}<+\infty$ hence $(\f_j)$ is uniformly
bounded from above. It follows from proposition 1.6 that $(\f_j)$ is relatively
compact. Indeed it can not converge uniformly to $-\infty$ since
$\sup_K \f_j=0$ (see proposition 1.6).
\end{proof}

\begin{pro} Let $K$ be a Borel subset of $X$.

1) If $K' \subset K$ then $V_{X,\o}  \leq V_{K,\o} \leq V_{K',\o}$
and $\sup_X V_{X,\o}=0$. Furthermore $V_{X,\o} \equiv 0$ when $\o \geq 0$.

2) If $\o_1 \leq \o_2$ then $V_{K,\o_1} \leq V_{K,\o_2}$.

3) For all $A>0$, $V_{K,A \o}=A \cdot V_{K,\o}$.

4) If $\o'=\o+dd^c \chi $ then
$$
-\chi+\inf_X \chi +V_{K,\o} \leq V_{K,\o'} \leq 
V_{K,\o}+\sup_X \chi -\chi.
$$

5) If $f:X \rightarrow X$ is holomorphic then
$$
V_{f(K),\o} \circ f \leq V_{K,f^*\o}.
$$
In particular if $f$ is a $\o$-isometry then $V_{f(K),\o}=V_{K,\o}$.
\end{pro}

\begin{proof}
That $K \mapsto V_{K,\o}$ is decreasing follows straightforwardly from the
definition. Observe that $0 \in PSH(X,\o)$ when $\o \geq 0$, hence
$V_{X,\o} \equiv 0$ in this case. When $\o$ is smooth (but not positive),
considering $V_{X,\o}^*$ will be a useful way of constructing
a positive closed current $\o_{V_{X,\o}^*} \sim \o$ with minimal
singularities (see section 5).

Assertions 2,3,4 are simple consequences of proposition 1.3.
The last assertion results from the following observation: if 
$\f \in PSH(X,\o)$ is such that $\f \leq 0$ on $f(K)$,
then $\f \circ f$ belongs to $PSH(X,f^* \o)$ and satisfies
$\f \circ f \leq 0$ on $K$.
\end{proof}

\begin{exa}
Assume $X=\C\P^n$, $\o$ is the Fubini-Study K\"ahler form and let
$B_R$ denote the euclidean ball centered at the origin and of radius $R$
in $\C^n \subset \C\P^n$. Then for $x \in \C^n$,
$$
V_{B_R,\o}(x)=\max\left( \log \frac{||x||}{R}+\frac{1}{2} \log[1+R^2]
-\frac{1}{2} \log[1+||x||^2];0 \right).
$$
Indeed set $\p_R:=\max(\frac{1}{2} \log[1+||x||^2],\f_R)$, where
$\f_R=\frac{1}{2} \log[1+R^2]+\log \frac{||x||}{R}$. Recall that the usual
Siciak's extremal function of $B_R$ is
$\log^+\frac{||x||}{R}$. Therefore
$\frac{1}{2} \log[1+||x||^2] \leq \f_R=\p_R$ for $||x|| \geq R$.
On the other hand if $||x||<R$ then 
$1+||x||^2>(1+R^2) \frac{||x||^2}{R^2}$ hence
$\frac{1}{2} \log[1+||x||^2]>\f_R$ in $B_R$.

Now let $u \in PSH(\C\P^n,\o)$ such that $u \leq 0$ in $B_R$.
Then $v=u+\frac{1}{2} \log[1+||x||^2] \in {\mathcal L}(\C^n)$. Since
$v \leq \frac{1}{2} \log[1+R^2]$ in $B_R$ we infer 
$v \leq \frac{1}{2} \log[1+R^2]+\log^+\frac{||x||}{R}=\p_R$ in 
$\C^n \setminus B_R$. Moreover $v \leq \frac{1}{2} \log[1+||x||^2]=\p_R$
in $B_R$ hence $v \leq \p_R$ in $\C^n$. This
shows $V_{B_R,\o}=\p_R-\frac{1}{2} \log[1+||x||^2]$ on $\C\P^n$.
\end{exa}

\begin{pro}
\text{ }

1) If $E$ is an open subset, then $V_E=V_E^*$.

2) Let $E$ be a Borel subset and $P$ a $PSH(X,\o)$-polar set. Then
$$
V_{E \cup P}^* \equiv V_E^*.
$$

3) Let $(E_j)$ be an increasing sequence of Borel subsets and 
set $E=\cup E_j$. Then
$V_{E,\o}^*=\lim \searrow V_{E_j,\o}^*$ if $\o$ is K\"ahler.

4) Let $K_j$ be a decreasing sequence of compact subsets of $X$ 
and set $K=\cap K_j$. Then
$V_{K_j,\o} \nearrow V_{K,\o}$, hence
$V_{K_j,\o}^* \nearrow V_{K,\o}^*$ a.e. 

5) Fix $E \subset X$ a non-pluripolar set. 
Then there exists $G_j$ a decreasing sequence of
open subsets, $E \subset G_j$, such that
$V_E^*=\lim V_{G_j}^*$.
\end{pro}

\begin{proof}
We write here $V_E$ for $V_{E,\o}$ since $\o$ is fixed and no
confusion can arise.

Let $E$ be an open subset of $X$. Observe that $V_E \leq 0$ on $E$,
hence $V_E^* \leq 0$ on $E$ which is open. Therefore $V_E^* \leq V_E$,
whence equality. This proves 1).

Let $w \in PSH(X,\o)$, $w \leq 0$, and fix $P \subset \{w=-\infty \}$.
Fix $E$ a Borel subset of $X$. Clearly $V_{E \cup P} \leq V_E$ hence
$V_{E \cup P}^* \leq V_E^*$. Conversely let $\f \in PSH(X,\o)$ be
such that $\f \leq 0$ on $E$. Then $\forall \e>0$,
$\p_{\e}:=(1-\e)\f+\e w \in PSH(X,\o)$ satisfies $\p_{\e} \leq 0$
on $E \cup P$, hence $\p_{\e} \leq V_{E \cup P}$. Letting
$\e \rightarrow 0$ we infer $\f \leq V_{E \cup P}$ on $X \setminus P$,
hence $\f \leq V_{E \cup P}^*$ on $X$. Thus $V_E^* \leq V_{E \cup P}^*$.

Let $E_j$ be an increasing sequence of subsets of $X$ and set 
$E= \cup_{j \geq 1} E_j$. Let $v:=\lim\searrow V_{E_j}^*$ (the limit is
decreasing by 4.4.1).  If $E$ is $PSH(X,\o)$-polar then so are all the
$E_j'$s, hence $V_E^* \equiv +\infty=\lim V_{E_j}^*$. So let us assume
$E$ is not $PSH(X,\o)$-polar. 
Then $v \in PSH(X,\o)$  since 
$v \geq V_{E,\o}^* \not\equiv-\infty$ (see proposition 1.6.3).
Observe that $v=0$ on the set $E \setminus N$,
where $N=\cup_{j \geq 1} \{V_{E_j}<V_{E_j}^*\}$. The latter is
called a {\it negligible set}. It follows from the local theory [5]
together with theorem 6.2 that $N$ is $PSH(X,\o)$-polar.
Therefore $V_E^* \leq v \leq V_{E \setminus N}^* =V_E^*$ by 2).

Let $K_j$ be a decreasing sequence of compact subsets and set
$K=\cap_j K_j$. Clearly $\lim \nearrow V_{K_j} \leq V_K$. Fix $\e>0$ and let
$\f \in PSH(X,\o)$ be such that $\f \leq 0$ on $K$.
Then $\{\f <\e\}$ is an open set which contains all $K_j's$, for
$j \geq j_{\e}$ large enough. Thus $\f-\e \leq 0$ on $K_j$, hence
$\f-\e \leq \lim \nearrow V_{K_j}$. Taking the supremum over all such $\f'$s
and 
letting $\e \rightarrow 0$ yields the reverse inequality
$V_K \leq \lim \nearrow V_{K_j}$. The conclusion on the convergence
of the upper semi-continuous regularizations follows now from
proposition 1.6.

It remains to prove 5). By Choquet's lemma, there exists an increasing
sequence $\f_j \in PSH(X,\o)$ such that $\f_j \leq 0$ on $E$ and
$V_E^*=(\sup_j \f_j)^*$. Set $G_j:=\{\f_j<1/j\}$. This defines a decreasing
sequence of open subsets containing $E$. Observe that
$\f_j-1/j \leq V_{G_j} \leq V_E$, hence
$\lim \f_j \leq \lim V_{G_j} \leq V_E$. Therefore
$V_E^*=\lim V_{G_j}^*$.
\end{proof}

\subsection{Alexander capacity}

\begin{defi}
Let $K$ be a Borel subset of $X$. We set
$$
T_{\o}(K):=\exp(-\sup_X V_{K,\o}^*).
$$
\end{defi}

This capacity characterizes again $PSH(X,\o)$-polar sets:

\begin{pro}
Let $P$ be a Borel subset. Then $T_{\o}(P)=0$ iff $P$ is $PSH(X,\o)$-polar.
Moreover if $\f \in PSH(X,\o)$ then
$$
T_{\o}(\f<-t) \leq C_{\f} \exp(-t), \; \forall t \in \R,
$$
where $C_{\f}=\exp(-\sup_X \f)$.
\end{pro}

\begin{proof}
The first assertion follows from theorem 4.2.
Let $\f \in PSH(X,\o)$, $t \in \R$ and set
$K_t=\{\f<-t\}$. Then $\f+t \leq 0$ on $K_t$ hence $\f+t \leq V_{K_t,\o}^*$.
We infer $\sup_X \f +t \leq \sup_X V_{K_t,\o}^*$ which yields
$T_{\o}(K_t) \leq \exp(-\sup_X \f) \exp(-t)$.
\end{proof}

The following proposition is an immediate consequence of proposition 4.4.
It shows that capacities $T_{\o},T_{\o'}$
are comparable if $\o,\o'$ are both K\"ahler. Further they enjoy nice 
invariance properties.

\begin{pro} \text{ }

1) For all Borel subsets $K' \subset K \subset X$, 
$T_{\o}(K') \leq T_{\o}(K) \leq T_{\o}(X)=1$.

2) If $\o_1 \leq \o_2$ then $T_{\o_1}(\cdot) \geq T_{\o_2}(\cdot)$.
Forall $A>0$, $T_{A \o}(\cdot)=[T_{\o}(\cdot)]^A$.
In particular if $\o$ and $\o'$ are both K\"ahler then
there exists $C \geq 1$ such that
$$
[T_{\o}(\cdot)]^C \leq T_{\o'}(\cdot) \leq [T_{\o}(\cdot)]^{1/C}.
$$

3) If $\o'=\o+dd^c \chi$ then 
$$
\frac{1}{C}T_{\o}(\cdot) \leq T_{\o'}(\cdot) \leq C \cdot T_{\o}(\cdot),
$$
where $C=\exp(\sup_X \chi-\inf_X \chi) \geq 1$.

4) If $f:X \rightarrow X$ is a holomorphic map then 
$T_{f^*\o}(\cdot) \leq T_{\o} \circ f(\cdot)$. In particular
if $f$ is a $\o$-isometry then $T_{\o} \circ f=T_{\o}$.
\end{pro}

\begin{rqe}
Following Zeriahi [43] one can prove that for all $\a < 2/\nu(X,\o)$
there exists $C_{\a}>0$ such that
$$
\text{Vol}_{\o}(\cdot) \leq C_{\a} T_{\o}(\cdot)^{\a},
$$
where $\nu(X,\o)=\sup \{ \nu(\f,x) \, / \, \f \in PSH(X,\o), x \in X \}$
and $\nu(\f,x)$ denotes the Lelong number of $\f$ at point $x$.
In particular it follows from proposition 4.8 that $\forall \f \in PSH(X,\o)$
with $\sup_X \f=0$,
$$
\text{Vol}_{\o}(\f<-t) \leq C_{\a} \exp(-\a t), \; \forall t \in \R.
$$
Such inequalities are quite useful in complex dynamics [22],[24] and in the
study of the complex Monge-Amp\`ere operator [30].
\end{rqe}

\begin{exa} 
Assume $X=\C\P^n$, $\o$ is the Fubini-Study K\"ahler form
and $B_R$ is the euclidean ball centered at the origin and of radius
$R$ in a chart $\C^n \subset \C\P^n$. We have explicitly computed 
the extremal function in this case (example 4.5). This yields
$$
T_{\o}(B_R)=\frac{R}{\sqrt{1+R^2}}.
$$
Observe that $T_{\o}(B_R) \sim R$ as $R \rightarrow 0$. This shows the
optimality of the rate of decreasing
in proposition 4.8.
\end{exa}

The capacity $T_{\o}$ in example 4.11
has to be related to the capacity $T_{\B^n}$
which measures compact subsets of the unit ball $\B^n$ of $\C^n$.
It is defined as follows: given $K$ a Borel subset of $\C^n$,
$T_{\B^n}(K):=\exp(-\sup_{\B^n} L_K)$, where
$$
L_K(z)=\sup \{ v(z) \, / \, v \in {\mathcal L}(\C^n), \; \sup_K v \leq 0 \}
$$
is the Siciak's extremal function of $K$ and ${\mathcal L}(\C^n)$ denotes
the Lelong class of psh functions with logarithmic growth in $\C^n$
(see example 1.2).
Let $\o=\o_{FS}$ denote the Fubini-Study K\"ahler form
on $\C\P^n$. One easily checks that
$$
V_{K,\o}-\log \sqrt{2} \leq L_K-\frac{1}{2}\log[1+|z|^2]
\leq V_{K,\o} \text{ in } \C^n. 
$$
We infer straightforwardly 
$\sup_{\B^n} L_K \leq \log \sqrt{2}+\sup_{\C\P^n}V_{K,\o}$ hence
$T_{\B^n}(K) \geq 2^{-1/2} T_{\o}(K)$. We also have a reverse inequality.
Indeed 
$\forall \f \in PSH(\C\P^n,\o)$, $\sup_{\C\P^n} \f \leq \sup_{\B^n} \f +C_1$,
where $C_1=\sup_{\C\P^n} V_{\B^n,\o}=\log \sqrt{2}$.
Therefore 
$$
\sup_{\C\P^n} V_{K,\o} \leq \sup_{\B^n} V_{K,\o} +\log \sqrt{2}
\leq \sup_{\B^n} L_K +\log 2,
$$
which yields
$$
\frac{1}{\sqrt{2}} T_{\o}(K) \leq
T_{\B^n}(K) \leq 2 T_{\o}(K).
$$

\begin{exa}
Assume again $X=\C\P^n$ and $\o$ is the Fubini-Study K\"ahler form.
Consider the totally real subspace $\R\P^n$ of points with real 
coordinates (the closure of $\R^n \subset \C^n$ in $\C\P^n$). Then
$$
\frac{1}{2(1+\sqrt{2})} \leq T_{\o}(\R\P^n) \leq 1.
$$
Indeed set $B_{\R^n}:=\R^n \cap \B^n$. 
It follows from the discussion above that 
$$
T_{\o}(\R\P^n) \geq \frac{1}{2}T_{\B^n}(\B_{\R^n}).
$$
Now there is an explicit formula
for $L_{B_{\R^n}}^*$ (Lundin's formula, see [29]),
$$
L_{B_{\R^n}}^*(z)=\sup\{ \log^+ |h(<z,\xi>)| \, / ||\xi||=1 \}, \; z \in \C^n,
$$
where $h(\zeta)=\zeta+\sqrt{\zeta^2-1}$. A simple computation yields
$|h(\zeta)| \leq \log[|z|+\sqrt{|z|^2+1}]$ for $\zeta=<z,\xi>$ with
$||\xi||=1$. We infer
$$
L_{\B_{\R^n}}(z) \leq \log\left[|z|+\sqrt{|z|^2+1}\right]
\leq \log[1+\sqrt{2}] \; \text{ in } \B^n,
$$
which yields the desired inequality.

Observe that the minorant is independent of the dimension $n$.
This has been used recently in complex dynamics by Dinh and Sibony [20].
\end{exa}

\begin{rqe}
It follows from proposition 4.6 that $T_{\o}$ is a generalized capacity
in the sense of Choquet which is outer regular.
\end{rqe}

\section{Tchebychev constants}

In this section we consider the case where $\o$ is (smooth
and) represents the first Chern class of a holomorphic line bundle
$L$ on $X$. 

Recall that a holomorphic line bundle $L$ on $X$ 
is a family of complex lines $\{L_x\}_{x \in X}$ together
with a structure of complex manifold of dimension $1+\dim_{\C} X$
such that the projection map $\pi:L \rightarrow X$ taking
$L_x$ on $x$ is holomorphic. Moreover one can always
locally trivialize $L$:
there exists an open covering $\{ {\mathcal U}_{\a} \}$
of $X$ and biholomorphisms
$\Phi_{\a}:\pi^{-1}({\mathcal U}_{\a}) \rightarrow {\mathcal U}_{\a} \times \C$
which take $L_x=\pi^{-1}(x)$ isomorphically onto $\{x\} \times \C$.
The line bundle $L$ is then uniquely 
(i.e. up to isomorphism) determined
by its transition functions 
$g_{\a \b} \in {\mathcal O}^*({\mathcal U}_{\a \b})$, 
${\mathcal U}_{\a \b}:={\mathcal U}_{\a} \cap {\mathcal U}_{\b}$,
where
$$
g_{\a\b}:=(\Phi_{\a} \circ \Phi_{\b}^{-1})_{|\{x\} \times \C}.
$$
Note that the $g_{\a\b}$'s satisfy the cocycle
condition $g_{\a \b} \cdot g_{\b \g} \cdot g_{\g \a} \equiv 1$,
hence define a class  $[\{g_{\a\b}\}] \in H^1(X,{\mathcal O}^*)$. 
The first Chern class of $L$ is the image $c_1(L) \in H^2(X,\Z)$
of $[\{g_{\a\b}\}]$ under the mapping 
$c_1:H^1(X,{\mathcal O}^*) \rightarrow H^2(X,\Z)$ induced
by the exponential short exact sequence 
$0 \rightarrow \Z \rightarrow {\mathcal O} \rightarrow {\mathcal O}^*
\rightarrow 0$. 

We let ${\bf \Gamma(X,L)}$ denote the set of holomorphic sections
of $L$ on $X$: $s \in \Gamma(X,L)$ is a collection $s=\{s_{\a}\}$
of holomorphic functions $s_{\a}$ on ${\mathcal U}_{\a}$ satisfying
the compatibility condition $s_{\a}=g_{\a\b} s_{\b}$ on 
${\mathcal U}_{\a\b}$. 
Similarly a (singular) metric $\p$ of $L$ on
$X$ is a collection $\p=\{\p_{\a}\}$ of functions
$\p_{\a} \in L^1({\mathcal U}_{\a})$ satisfying 
$\p_{\a}=\p_{\b}+\log|g_{\a\b}|$ in ${\mathcal U}_{\a\b}$.
The metric is said to be {\it smooth}
if the $\p_{\a}'s$ are ${\mathcal C}^{\infty}$-smooth functions.
A smooth metric always exists.
The metric $\p$ is said to be {\it positive} if the $\p_{\a}$'s
are psh functions. In particular if $s=\{s_{\a} \}$
is a holomorphic section of $L$ on $X$, then 
$\p=\{ \p_{\a}:=\log|s_{\a}| \}$ is a positive (singular) metric
of $L$ on $X$. Note that we make here a slight abuse of terminology:
differential geometers usually call "metric"  the
non-negative (usually smooth and non vanishing) quantities
$e^{-\p}=\{ e^{-\p_{\a}} \}$.

Given a (singular) metric $\p=\{\p_{\a}\}$ of $L$ on $X$, we consider
its curvature $\Theta_{\p}:=dd^c \p_{\a}$ in ${\mathcal U}_{\a}$.
This yields a globally well defined real closed current
on $X$ since $dd^c \log|g_{\a \b}|=0$ in ${\mathcal U}_{\a\b}$.
It is a standard consequence of de Rham's isomorphism 
that this current represents the image of the first Chern class
of $L$ under the mapping $i:H^2(X,\Z) \rightarrow H^2(X,\R)$
(induced by the inclusion $i:\Z \rightarrow \R$).
The line bundle $L$ is said to be {\it pseudoeffective} 
(resp. {\it positive}) if it admits a (singular) positive metric
(resp. a smooth metric whose curvature is a K\"ahler form).

Fix $h=\{h_{\a}\}$ a smooth metric of $L$ on $X$ and set
$\o:=\Theta_h$. Then $PSH(X,\o)$
is in $1$-to-$1$ correspondence with the set of positive
singular metrics of $L$ on $X$. Indeed if $\p$ is such a metric
then $\f:=\p-h$ is globally well defined on $X$ and such that
$dd^c \f \geq -\o$.
Conversely if $\f \in PSH(X,\o)$ then $\p=\{\p_{\a}:=\f+h_{\a}\}$
defines a positive singular metric of $L$ on $X$.
We can thus rephrase the pseudoeffectivity property as follows:
$$
L \text{ is pseudoeffective } \Longleftrightarrow PSH(X,\o) \neq \emptyset.
$$

Given $L$ a pseudoeffective line bundle, it is interesting
to know whether $L$ admits a positive metric which is less singular
than any another. This notion has been
introduced in [18] and happens to be related to very special
extremal functions:

\begin{pro}
Let $(L,h)$ be a pseudoeffective line bundle on $X$ equipped with
a smooth metric $h$. Set $\o:=\Theta_h$. Then
$$
h_{min}:=h+V_{X,\o}^*
$$
is a positive singular metric of $L$ on $X$ with "minimal singularities".
More precisely if $\p$ is a positive singular metric of $L$ on $X$,
then there exists a constant $C_{\p}$ such that
$\p \leq h_{min}+C_{\p}$.
\end{pro}

\begin{proof}
Let $\p$ be a positive singular metric of $L$ on $X$. Then
$\p-h$ is a globally well defined $\o$-psh function. It is 
u.s.c. hence bounded from above on $X$: we let $C_{\p}$ denotes its 
maximum. Then $\p-h-C_{\p} \leq 0$ on $X$, hence
$\p-h \leq V_{X,\o}^*+C_{\p}$, which yields $\p \leq h_{min}+C_{\p}$. 
\end{proof}

In the sequel we assume $L$ is positive and $h$ has been chosen
so that $\o:=\Theta_h$ is a K\"ahler form.
For $s \in \Gamma(X,L^N)$, we let $||s||_{Nh}$ denote the norm of $s$
computed with respect to the metric $N h$: it is defined
in ${\mathcal U}_{\a}$ by
$||s||_{Nh}:=|s_{\a}|e^{-N h_{\a}}$. The definition is independent
of $\a$ thanks to the compatibility conditions.

For a given Borel subset $K$ of $X$, we define its {\it Tchebychev constants}
$$
M_{N \o}(K):=\inf \left\{ \sup_K ||s||_{N h} \, / \, 
s \in \Gamma(X,L^N), \, \sup_X ||s||_{N h}=1 \right\}.
$$
Note that an obvious rescaling argument shows that 
$M_{d \o}$ remains unchanged if we replace $h$ by $h+C$
so that it really depends on $\o=\Theta_h$ rather than on $h$.
Consider
$$
T_{\o}'(K):=\inf_{N \geq 1} [M_{N \o}(K)]^{1/N}.
$$

\begin{thm}
Let $K$ be a compact subset of $X$. Then
$$
T_{\o}(K)=T_{\o}'(K).
$$
\end{thm}

\begin{proof}
The core of the proof consists in showing that
$$
V_{K,\o}(x)=\sup \left\{ \frac{1}{N} \log ||s||_{Nh}(x) \, / \, N \geq 1,
s \in \Gamma(X,L^N) \text{ and } \sup_K ||s||_{Nh} \leq 1 \right\}.
$$
Note that for any of the sections $s$ involved in the supremum,
$\f:=N^{-1}\log ||s||_{Nh}$ belongs to $PSH(X,\o)$ and satisfies
$\f \leq 0$ on $K$.
Therefore $\f \leq V_{K,\o}$.

Conversely fix $x_0 \in X$ and $a<V_{K,\o}(x_0)$. Fix $\f \in PSH(X,\o)$
such that $\sup_K \f \leq 0$ and $\f(x_0)>a$. Regularizing $\f$
(see Appendix) and translating, we can assume
$\f \in PSH(X,\o) \cap {\mathcal C}^{\infty}(X)$, $\sup_K \f<0$ and
$\f(x_0)>a$. 
Fix $\e>0$. Let $B=B(x_0,r)$ be a small ball on which $\f>a$. 
We choose $B$ so small that the oscillation of $h$ is smaller
than $\e$ on $B$.
Let $\chi$ be a test function with compact support in $B$
and such that $\chi \equiv 1$ in $B(x_0,r/2)$. 
We can assume w.l.o.g. that $B \subset {\mathcal U}_{\a_0}$
for some $\a_0$ but $B \cap {\mathcal U}_{\b}=\emptyset$
for all $\b \neq \a_0$. This insures that $\chi$ is a smooth
section of $L^N$ for all $N \geq 1$.

Let $\p_1$ be  a smooth
positive metric of $L^{N_1} \otimes K_X^*$ on $X$ (this is possible if
$N_1$ is chosen large enough since $L$ is positive). 
Let $\p_2$ be a positive metric of $L^{N_2}$ on $X$ which is smooth in
$X \setminus \{x_0\}$ and with Lelong number 
$\nu(\p_2,x_0) \geq n=\dim_{\C} X$ (this is again possible if $N_2$ is large
enough, since $L$ is ample). Observe that $\D \chi$ is a smooth
$\D$-closed $(0,1)$-form with values in $L^N$ (for all $N \geq 1$).
Alternatively it is a smooth $\D$-closed $(n,1)$-form with values
in $L^N \otimes K_X^*$. Applying H\"ormander's $L^2$-estimates
(see e.g. [15], chapter VIII)
with weight $\p_N:=(N-N_1-N_2)(\f+h)+\p_1+\p_2$, we find a smooth section
$f$ of $L^N$ such that $\D f= \D \chi$ and
$$
\int_X |f|^2 e^{-2(N-N_1-N_2)(\f+h)-2\p_1-2\p_2}dV_{\o}
\leq C_1 \int_X |\D \chi|^2 e^{-2\p_N} dV_{\o}.
$$
Note that $\D \chi$ has support in $B \setminus B(x_0,r/2)$ where
$\p_N$ is smooth so that both integrals are finite.
Since $\nu(\p_2,x_0) \geq n$, this forces $f(x_0)=0$. The second
integral is actually bounded from above by $C_2 e^{-2N(a-\e)}$,
where $C_2$ is independent of $N$, since $-\f<-a$ on $B$
and the oscillation of $h$ is smaller than $\e$ on $B$.
Therefore $s:=\chi-f \in \Gamma(X,L^N)$ satisfies $s(x_0)=1$
and 
$$
\int_X |s|^2e^{-2N(\f+h)} dV_{\o} \leq C_3 e^{-2N(a-\e)},
$$
where $C_3$ is independent of $N$. Now $\f <0$ in a neighborhood of $K$, so
the mean-value inequality applied to the subharmonic functions
$|s_{\a}|^2$ yields for all $x$ in $K$,
\begin{eqnarray*}
|s|^2e^{-2Nh}(x) &\leq& C_{\d} \int_{B(x,\d)} |s|^2(y)
e^{-2N[\f+h](y)} e^{2N[h(y)-h(x)+\f(y)]}d\l(y) \\
&\leq& C_4 e^{-2N(a-\e)}
\end{eqnarray*}
if $\d$ is so small that $|\sup_{B(x,\d)} \f|>0$ is bigger
than the oscillation of $h$ on $B(x,\d)$.
Therefore $S:=C_4^{-1/2} e^{N(a-\e)}s \in \Gamma(X,L^N)$ satisfies
$\sup_K ||S||_{Nh}  \leq 1$ and 
$N^{-1} \log ||S||_{Nh}(x_0) \geq a-\e-\frac{\log C_4}{2N}$.
Letting $N \rightarrow +\infty$, $\e \rightarrow 0$ and 
$a \rightarrow V_{K,\o}(x_0)$ completes the proof of the equality.

To conclude observe that by rescaling one gets
\begin{eqnarray*}
\lefteqn{-\log T_{\o}(K)=\sup_X V_{K,\o}} \\
&=& \! \! \! \!
\sup \left\{ \frac{1}{N} \sup_X \log ||S||_{N h} \, / \, N \geq 1, 
S \in \Gamma(X,L^N) \text{ and } \sup_K ||S||_{N h}=1 \right\} \\
&=&\! \! \! \!
\sup \left\{ -\frac{1}{N} \sup_K \log ||S||_{N h} \, / \, N \geq 1, 
S \in \Gamma(X,L^N) \text{ and } \sup_X ||S||_{N h}=1 \right\} \\
&=&-\log T_{\o}'(K).
\end{eqnarray*}
\end{proof}

\noindent {\bf Projective capacity.}
We assume here that $X=\C\P^n$ is the complex projective space
and $\o=\o_{FS}$ is the Fubini-Study K\"ahler form. We give in this
context a geometrical interpretation of the capacity
$T_{\o}$. This will shed some light on the notion of projective
capacity introduced by Alexander [1].

Let $\pi:\C^{n+1} \setminus \{ 0 \} \rightarrow \C\P^n$ denote the canonical 
projection map. We let $\B^{n+1}$ denote the unit ball in $\C^{n+1}$.
Recall that the polynomially convex hull $\widehat{F}$ of a compact
set $F$ of $\C^{n+1}$ is defined as 
$\widehat{F}:=\{x \in \C^{n+1} \, / \, |P(x)| \leq \sup_F |P|, \, \forall P \text
{ polynomial}\}$.

The following result gives an interesting interpretation of the
capacity $T_{\o}$.

\begin{thm}
Let $K$ be a compact subset of $\C\P^n$. Then
$$
T_{\o}(K)=\sup \{ r>0 \, / \, r \B^{n+1} \subset \widehat{K_0} \},
$$
where $K_0=\pi^{-1}(K) \cap \partial \B^{n+1}$.
\end{thm}

\begin{proof}
Let $K,K_0$ be as in the theorem. Observe that $K_0$ is a circled subset
of $\partial \B^{n+1}$: if $z \in K_0$ then $e^{i\theta}z \in K_0$,
$\forall \theta \in [0,2\pi]$. For such compacts, the polynomial hull
$\widehat{K_0}$ coincides with the "homogeneous polynomial hull",
$$
\widehat{K_0}^h:=\{x \in \C^{n+1} \, / \, |P(x)| \leq \sup_F |P|, \, 
\forall P \text{ homogeneous polynomial}\}.
$$
Indeed one inclusion $\widehat{K_0} \subset \widehat{K_0}^h$ is clear, so
assume $z_0 \in \widehat{K_0}^h$. Let $P=\sum_{j=0}^d P_j$ be a polynomial
of degree $d$ decomposed into its homogenous components. Observe that
$P_j(x)=(2 \pi)^{-1} \int_0^{2 \pi} P(e^{i\theta}x)e^{-ij\theta}d\theta$.
Therefore $\sup_{K_0} |P_j| \leq \sup_{K_0} |P|$ since $K_0$ is circled.
Fix $t \in ]0,1[$. Then
$$
|P(tz_0)| \leq \sum_{j=0}^d t^j |P_j(z_0)| \leq \frac{1}{1-t} 
\sup_{K_0} |P|.
$$
We infer $t z_0 \in \widehat{K_0}$. Letting $t \rightarrow 1^-$ and using that
$K_0$ is closed we get $z_0 \in \widehat{K_0}$, whence $\widehat{K_0}= \widehat{K_0}^h$.

Fix now $z \in \C^{n+1}$ such that $||z|| \leq T_{\o}(K)$. Let $P$ be
a homogeneous polynomial of degree $d$. Then
\begin{equation}
|P(z)|=||z||^d \left|P\left(\frac{z}{||z||}\right)\right|
\leq T_{\o}(K)^d \sup_{\partial \B^{n+1}} |P| 
\end{equation}
Now set $\p(z)=d^{-1} \log|P(z)|-\log||z||$ and $\f=\p-\sup_K \p$.
Then $\f \in PSH(X,\o)$ with $\sup_K \f \leq 0$ hence
$\f \leq V_{K,\o}$. Therefore
$$
T_{\o}(K)^d \leq \exp(-d \sup_{\C\P^n} \f)=
\frac{\sup_{K_0} |P|}{\sup_{\partial \B^{n+1}} |P|}.
$$
Together with $(2)$ this yields $|P(z)| \leq \sup_{K_0} |P|$
hence $z \in \widehat{K_0}^h=\widehat{K_0}$. Thus $K_0$ contains the ball
centered at the origin of radius $T_{\o}(K)$.

Conversely since $T_{\o}(K)=T_{\o}'(K)$ (theorem 4.1), one can find
homogenous polynomials $P_j$ of degree $d_j$ such that
$\sup_{\partial \B^{n+1}} |P_j|^{-1/d_j} \cdot \sup_{K_0} |P_j|^{1/d_j}
\rightarrow T_{\o}(K).$ Assume $r \B^{n+1} \subset \widehat{K_0}$. Then
$$
r^{d_j} \sup_{\partial \B^{n+1}} |P_j|=\sup_{r \B^{n+1}} |P_j| 
\leq \sup_{K_0} |P_j|
$$
yields $r \leq T_{\o}(K)$.
\end{proof}

\begin{rqe}
Sibony and Wong [36]
have been first in showing that if a compact subset $K$ of
$\C\P^n$ is large enough then the polynomial hull of $K_0$
contains a full neighborhood of the origin in $\C^{n+1}$.
They used the (complicated) notion of $\Gamma$-capacity. Their
approach has been simplified by Alexander [1] who
introduced a projective capacity which is comparable
to $T_{\o}$ (see theorem 4.4 in [1]). The proof given above
is essentially Alexander's (see also theorem 4.3 in [38]).

This result has been used recently in complex
dynamics (see [19],[25]).
\end{rqe}

\noindent {\bf Further capacities.}
In our definition of Chebyshev constants we have normalized
holomorphic sections $s \in \Gamma(X,L^N)$ by requiring
$\sup_X ||s||_{N h}=1$. Given $\mu$ a probability measure
such that $PSH(X,\o) \subset L^1(\mu)$ and $A \in \R$, we
could as well consider
$$
M_{N,\o}^{\mu,A}(K):=\left\{ \sup_K ||s||_{N h} \, / \, 
s \in \Gamma(X,L^N), \, \int_X \log ||s||_{Nh} d\mu=A \right\}.
$$
This normalization has the following pleasant property: if
$s \in \Gamma(X,L^N)$ and $s' \in \Gamma(X,L^{N'})$ are so 
normalized then $s \cdot s' \in \Gamma(X,L^{N+N'})$ again
satisfies $\int_X \log ||ss'||_{(N+N')h} d\mu=A$. 
We infer
$M_{N+N',\o}^{\mu,A} \leq M_{N,\o}^{\mu,A} \cdot M_{N',\o}^{\mu,A}$
so that
$$
T_{\o}^{\mu,A}(K):=\inf_{N \geq 1}[M_{N,\o}^{\mu,A}(K)]^{1/N}
=\lim_{N \rightarrow +\infty} [M_{N,\o}^{\mu,A}(K)]^{1/N}.
$$
This yields a whole family of capacities which are
all comparable to $T_{\o}$ thanks to proposition 1.7: there
exists $C=C(\mu,A) \geq 1$ such that
$$
\frac{1}{C} T_{\o}(\cdot) \leq T_{\o}^{\mu,A}(\cdot)
\leq C T_{\o}(\cdot).
$$
The {\it projective capacity} of Alexander [1] is
precisely $T_{\o}^{\mu,A}$ for $X=\C\P^n$, $\o=\o_{FS}$,
$\mu=\o^n$ and 
$A=\int_{\C\P^n} \left(\log|z_n|-\log||(z_0,\ldots,z_n)|| \right) \o^n([z]).$

\section{Comparison of capacities and applications}

\subsection{Josefson's theorem}

In this section we assume that $\o$ is K\"ahler and normalized
by $Vol_{\o}(X)=1$. We first prove inequalities relating
$T_{\o}$ and $Cap_{\o}$. 
Then we prove (theorem 6.2) a quantitative version of Josefson's 
theorem that every locally pluripolar set is actually
$PSH(X,\o)$-polar. In the local theory this result is due to
El Mir [21]. We follow the approach of Alexander-Taylor [2].

\begin{pro}
There exists $A>0$ s.t. for all
compact subsets $K$ of $X$,
$$
\exp\left[-\frac{A}{Cap_{\o}(K)} \right]
 \leq T_{\o}(K) \leq e \cdot 
\exp\left[-\frac{1}{Cap_{\o}(K)^{1/n}}\right].
$$
\end{pro}

\begin{proof}
Set $M_K=\sup_X V_{K,\o}$. If $M_K=+\infty$ then $K$ is $PSH(X,\o)$-polar
(theorem 4.2) and there is nothing to prove:
$T_{\o}(K)=Cap_{\o}(K)=0$.
So we assume in the sequel $M_K<+\infty$ hence $V_{K,\o}^* \in PSH(X,\o)$.
If $M_K \geq 1$ then $u_K:=M_K^{-1} V_{K,\o}^* \in PSH(X,\o)$
with $0 \leq u_K \leq 1$ on $X$. Since $\o_{V_{K,\o}^*} \leq M_K \o_{u_K}$,
we get
$$ 
\frac{1}{M_K^n}=\frac{1}{M_K^n} \int_K (\o_{V_{K,\o}^*})^n 
\leq \int_K (\o_{u_K})^n \leq Cap_{\o}(K)
$$
whence $T_{\o}(K) \leq \exp(-Cap_{\o}(K)^{-1/n})$.

If $0 \leq M_K \leq 1$ then $0 \leq V_{K,\o}^* \leq 1$ hence
$V_{K,\o}-1$ coincides with the relative extremal function $h_{K,\o}$
(see proposition 2.14).
We infer
$$
1=\int_K (\o_{V_{K,\o}^*})^n \leq Cap_{\o}(K) \leq Cap_{\o}(X)=1,
$$
while $T_{\o}(K) \leq T_{\o}(X)=1$.
Thus in both cases 
$T_{\o}(K) \leq e \cdot \exp( Cap_{\o}(K)^{-1/n})$.

\vskip.2cm
We now prove the reverse inequality. We can assume $M_K \geq 1$,
otherwise it is sufficient to adjust the value of $A$.
Let $\f \in PSH(X,\o)$ be such that $\f \leq 0$ on $K$.
Then $\f \leq M_K$ on $X$, hence
$w:=M_K^{-1}(\f-M_K) \in PSH(X,\o)$ satisfies $\sup_X w \leq 0$ and 
$w \leq -1$ on $K$. We infer
$w \leq h_{K,\o}^*$, hence
$$
w_K:=\frac{V_{K,\o}^*-M_K}{M_K} \leq h_{K,\o}^* \leq 0.
$$
Now $\sup_X (V_{K,\o}^*-M_K)=0$, so it follows from proposition 1.7 that
$\int_X |V_{K,\o}^*-M_K| \o^n \leq C_1$ for some constant $C_1>0$ independent
of $K$. We infer
\begin{eqnarray*}
Cap_{\o}(K)=\int_K (\o_{ h_{K,\o}^*})^n &\leq&
\int_X [- h_{K,\o}^*] (\o_{ h_{K,\o}^*})^n \\
&\leq &\frac{1}{M_K} \int_X -(V_{K,\o}^*-M_K) (\o_{ h_{K,\o}^*})^n
\leq \frac{C_2}{M_K},
\end{eqnarray*}
using corollary 2.3 and the fact
that $ h_{K,\o}^*=-1$ on $K$, except perhaps on a pluripolar set which
has zero $(\o_{ h_{K,\o}^*})^n$-measure.
This yields the desired inequality.
\end{proof}

It follows from the previous proposition and corollary 2.8 that
$\o$-psh functions are quasicontinuous with respect to the capacity
$T_{\o}$.

\vskip.2cm

\begin{thm}
Locally pluripolar sets are 
$PSH(X,\o)$-polar. 
\end{thm}

\begin{proof}
More precisely we are going to show the following:
consider $\Omega$  an open subset of $X$,
$v \in PSH^-(\Omega)$ and $P \subset \{v=-\infty\}$.
Fix $0<\e<1/n$ and $V_t:=V_{G_t,\o}$ where
$G_t=\{ x \in \Omega \, / \, v(x)<-t \}$. Then
$$
\f_{\e}(x):=\frac{1}{\e} \int_1^{+\infty} \frac{1}{t^{1+\e}}
[V_t(x)-\sup_X V_t] dt
$$
is a $\o$-psh function such that $P \subset \{\f_{\e}=-\infty\}$.

Indeed since $G_t$ is open, 
we have $V_t \in PSH(X,\o)$ and $V_t=0$ on $G_t$
(see proposition 4.6).
Observe that $\f_{\e}$ is a sum of negative $\o$-psh functions
hence it is either identically $-\infty$ or a well defined $A \o$-psh
function with $A=\e^{-1} \int_1^{+\infty} t^{-(1+\e)}dt=1$.
Recall that $-C+\sup_X V_t \leq \int_X V_t \o^n \leq \sup_X V_t$
(proposition 1.7). Therefore $\int_X \f_{\e} \o^n \geq -C$
hence $\f_{\e} \in PSH(X,\o)$.

Fix $x \in \Omega$ such that $v(x)<-1$.
Observe that $V_t-\sup_X V_t \leq 0$ with $V_t(x)=0$ if
$x \in  G_t$, i.e. when $|v(x)|>t$. Therefore
$$
\f_{\e}(x) \leq -\frac{1}{\e} \int_1^{|v(x)|} \frac{\sup_X V_t}{t^{1+\e}}dt.
$$
Recall now that $Cap_{\o}$ is always dominated by $Cap_{BT}$ hence
$Cap_{\o}(G_t)\leq C_1/t <1$ if $t$ is large enough. We infer
from the previous proposition that
$$
-\sup_X V_t \leq -[Cap_{\o}(G_t)]^{-1/n} \leq -C_2 t^{1/n},
$$
which yields 
$$
\f_{\e}(x) \leq \frac{-C_2}{\e} \int_1^{|v(x)|} \frac{dt}{t^{1+\e-1/n}}
\leq -C_3 |v(x)|^{1/n-\e}+C_4.
$$
Note that $\f_{\e}(x)=-\infty$
whenever $v(x)=-\infty$ hence $P \subset \{\f_{\e}=-\infty \}$.
\end{proof}

\subsection{Dynamical capacity estimates}

Let $f:\C\P^n \rightarrow\C\P^n$ be an holomorphic endomorphism.
We let $\o$ denote again the Fubini-Study K\"ahler form.
Then $f^* \o$ is a smooth positive closed $(1,1)$-form of
mass $\l=\int_{\C\P^n} f^* \o \wedge \o^{n-1}=:$the 
{\it first algebraic degree of $f$}. 
Thus $\l^{-1}f^*\o=\o+dd^c \f$, where $\f$ is a smooth $\o$-psh 
function on $\C\P^n$.
Iterating this functional equation yields
$$
\frac{1}{\l^j}(f^j)^*\o=\o+dd^c g_j, \; 
g_j=\sum_{l=0}^{j-1} \frac{1}{\l^l} \f \circ f^l.
$$
We assume $\l \geq 2$. Thus the sequence $(g_j)$ uniformly converges
on $\C\P^n$ towards a continuous function $g_f \in PSH(X,\o)$
called the {\it Green function} of $f$. We refer the interested
reader to [35] for a detailed study of the properties of
the {\it Green current} $T_f=\o+dd^c g_f$.

Dynamical volume estimates have revealed quite useful in establishing ergodic
properties of the Green current $T_f$ 
(see [22], [24] and references therein). We 
establish herebelow very simple dynamical capacity estimates
and show how to derive from them dynamical volume estimates.

\begin{pro}
There exists $0<\a<1$ such that for all Borel subsets $K$ of $X$,
for all $j \in \N$,
$$
\left[ \a T_{\o}(K)\right]^{\l^j} \leq T_{\o}(f^j(K)).
$$
\end{pro}

\begin{proof}
This follows straightforwardly from proposition 4.9:
$$
T_{\o}(f^j(K)) \geq T_{(f^j)^*\o}(K)=
\left[T_{\l^{-j}(f^j)^*\o}(K)\right]^{\l^j} \geq 
\left[ \a T_{\o}(K)\right]^{\l^j},
$$
where the first two inequalities follow from 4.9.4 and 4.9.2 and
last one follows from 4.9.3 and the fact that
$\l^{-j}(f^j)^*\o=\o+dd^c g_j$, where $g_j$ is uniformly bounded.
\end{proof}

\begin{cor}
Let $\f \in PSH(X,\o)$. Then the sequence
$(\l^{-j} \f \circ f^j)$ is relatively compact in $L^1(\C\P^n)$.
\end{cor}

\begin{proof}
Set $\f_j=\l^{-j} \f \circ f^j$. Observe that $\f_j$ is uniformly bounded
from above and that $\f_j+g_j \in PSH(\C\P^n,\o)$.
It follows from proposition 1.6 that either
$\f_j$ converges uniformly towards $-\infty$ or it is relatively compact
in $L^1(\C\P^n)$.
It is sufficient to show that for $A>0$ large enough,
$\overline{\lim}_{j \rightarrow +\infty} T_{\o}(\f_j<-A)<T_{\o}(X)=1$.
Observe that $f^j(\f_j<-A)=\{\f<-A\l^j\}$. Therefore
$$
\left[ \a T_{\o}(\f_j<-A)\right]^{\l^j} \leq T_{\o}(\f<-A\l^j)
\leq C \exp(-A \l^j),
$$
where the last inequality follows from proposition 4.8. We infer
$$
\overline{\lim}_{j \rightarrow +\infty} T_{\o}(\f_j<-A)
\leq \frac{1}{\a} \exp(-A)<1
$$
for $A>-\log \a$ large enough.
\end{proof}

\begin{cor}
There exists $C>0$ such that for all Borel subset
$K$ of $\C\P^n$ and for all $j \in \N$, 
$$
Vol_{\o}(f^j(K)) \geq \exp \left( -\frac{C \l^j}{Vol_{\o}(K)} \right).
$$
\end{cor}

In other words the volume of a given set can not decrease too fast
under iteration. Such volume estimates are used in complex dynamics
to prove fine convergence results towards the Green current $T_f$ (see [22],
[24]). One may hope that dynamical capacity estimates will allow
to establish convergence results in higher codimension.

\begin{proof}
By the change of variables formula one gets
$$
Vol_{\o}(f^jK)=\int_{f^jK} \o^n \geq 
\frac{1}{d_t^j} \int_K (f^j)^* \o^n=
\frac{1}{d_t^j} \int_K |J_{FS}(f^j)|^2 \o^n,
$$
where $d_t=\l^n$ denotes the topological degree of $f$ and 
$J_{FS}(f)$ stands for the jacobian of $f$ with respect to the Fubini-Study
volume form. Observe that
$\log |J_{FS}(f)|=u-v$ is a difference of two qpsh functions
$u,v \in PSH(X,A \o)$ for some $A=A(\l,f)$. 
Moreover by the chain rule,
$$
\frac{1}{\l^j}\log|J_{FS}(f^j)|=\sum_{l=0}^{j-1} \frac{1}{\l^j}
\log |J_{FS}(f) \circ f^l|.
$$
Since $\l^{-l} \log |J_{FS}(f) \circ f^l|$ is relatively compact
in $L^1(\C\P^n)$ (previous corollary), the concavity of the log
yields 
\begin{eqnarray*}
\frac{1}{Vol_{\o}(K)} \int_K |J_{FS}(f^j)|^2 \o^n &\geq&
\exp \left( \frac{2 \l^j}{Vol_{\o}(K)} 
\int_K \frac{1}{\l^j}\log|J_{FS}(f^j)| \o^n \right) \\
&\geq& \exp\left( -\frac{C_1 \l^j}{Vol_{\o}(K)} \right).
\end{eqnarray*}
The conclusion follows by observing that $\a \exp(-x/\a) \geq \exp(-2x/\a)$,
for all $\a>0$ and all $x \geq 1/e$.
\end{proof}

\section{Appendix: Regularization of qpsh functions}

It is well-known that every psh 
function $\f$ can be 
{\it locally regularized}, i.e. one can find {\it locally} a sequence $\f_j$
of smooth psh functions which decrease towards $\f$
(see e.g. [17], chapter 1). Similarly one can always 
locally regularize $\o$-psh functions. It is interesting
to know whether one can also {\it globally regularize}
$\o$-psh functions.

When $X$ is a complex homogeneous manifold (i.e. when $Aut(X)$ acts
transitively on $X$), it is possible to approximate any $\o$-psh
function by a decreasing sequence of smooth $\o$-psh functions
(see [23], [27]). In general however there is
a loss of positivity: it will be possible to approximate $\f \in PSH(X,\o)$
by a decreasing sequence of smooth functions $\f_j$ but
the curvature forms $dd^c \f_j$ will have to be more negative than $-\o$.
How negative depends on the positivity of the cohomology class
$[\o]$.

Consider e.g. $\pi:X \rightarrow \P^2$ the blow up of $\P^2$ at point
$p$, $E=\pi^{-1}(p)$ the exceptional divisor and let $\o=[E]$ be the
current of integration along $E$. Then $PSH(X,\o) \simeq \R$ 
(see Remark 1.5) so every psh function has logarithmic singularities
along $E$, hence is not smooth. Alternatively $E$ has self-intersection
$-1$ so its cohomology class cannot be represented by smooth non-negative
forms, not even by smooth forms with (very) small negativity.

Following Demailly's fundamental work [12], [14], [18]
(to cite a few)
we show herebelow that regularization with no loss of positivity
is possible when $\o$ is a Hodge form (i.e. a K\"ahler form
with integer class). This yields a "simple" regularization process
when $X$ is projective. We would like to mention that Demailly has
produced over the last twenty years much finer regularization results.
We nevertheless think it is worth including a proof, since it is far
less technical than Demailly's more general results
(although our proof heavily relies
on his ideas). We thank P.Eyssidieux for his helpful contribution
regarding that matter.

\begin{thm}
Let $L \rightarrow X$ be a positive holomorphic line bundle equipped
with a smooth strictly positive metric $h$, and set $\o:=\Theta_h >0$.

Then for every $\f \in PSH(X,\o)$, there exists a  
sequence $\f_j \in PSH(X,\o) \cap {\mathcal C}^{\infty}(X)$ such that
$\f_j$ decreases towards $\f$.
\end{thm}

\begin{proof}
Let $\f \in PSH(X,\o)$. We can assume w.l.o.g. that $\f \leq 0$ on $X$.
Let $\p=\{ \p_{\a}:=\f+h_{\a} \in PSH({\mathcal U}_{\a})\}$ 
denote the associated (singular)
positive metric of $L$ on $X$, where $\{{\mathcal U}_{\a}\}$ denotes an open
cover of $X$ trivializing $L$ (see section 4). 

{\bf Step 1.} We consider the following Bergman spaces
$$
{\mathcal H}_{j,j_0}:=\left\{ s \in \Gamma(X,L^j) \, / \, 
\int_X |s|^2 e^{-2h_{j,j_0}} dV_{\o}<+\infty \right\},
$$
where $h_{j,j_0}=(j-j_0)\p+j_0h$, $j_0$ a fixed large integer
(to be specified later). Let 
$\sigma_1^{(j,j_0)},\ldots,\sigma_{s_j}^{(j,j_0)}$ be an orthonormal basis
of ${\mathcal H}_{j,j_0}$ and set
$$
\p_{j,j_0}:=\frac{1}{2j} \log \left[ \sum_{l=1}^{s_j} |\sigma_l^{(j,j_0)}|^2
\right]=\frac{1}{2j} \sup_{s \in B_{j,j_0}} \log |s|^2,
$$
where $B_{j,j_0}$ denotes the unit ball of radius 1 centered at 0 in
${\mathcal H}_{j,j_0}$. Clearly $\p_{j,j_0}$ defines a positive (singular)
metric of $L$ on $X$, equivalently $\f_{j,j_0}:=\p_{j,j_0}-h \in PSH(X,\o)$.
If $x \in {\mathcal U}_{\a}$ and $s=\{s_{\a}\} \in {\mathcal H}_{j,j_0}$,
then $|s_{\a}|^2$ is subharmonic in ${\mathcal U}_{\a}$ hence
$$
|s_{\a}(x)|^2 \leq \frac{C_1}{r^{2n}} \int_{B(x,r)} |s_{\a}(x)|^2
\leq \frac{C_2}{r^{2n}} e^{2 \sup_{B(x,r)} h_{j,j_0}} 
\int_X |s|^2 e^{-2h_{j,j_0}} dV_{\o},
$$
where $r>0$ is  so small that $B(x,r) \subset {\mathcal U}_{\a}$.
We infer
\begin{equation}
\f_{j,j_0}(x) \leq (1-j_0/j) \sup_{B(x,r)} \f +\frac{C_3-n \log r}{j}.
\end{equation}

There is also a reverse inequality which uses a deep extension
result of Ohsawa-Takegoshi-Manivel (see [16]): there exists
$j_0 \in \N$ and $C_4>0$ large enough so that 
$\forall x \in X, \forall j \in \N$, there exists $s \in \Gamma(X,L^j)$
with
$$
\int_X |s|^2 e^{-2h_{j,j_0}} dV_{\o} \leq C_4 |s(x)|^2 e^{-2h_{j,j_0}(x)}.
$$
Choose $s$ so that the right hand side is equal to $1$,
hence $s \in B_{j,j_0}$. Then
$$
\p_{j,j_0}(x) \geq \frac{1}{2j} \log |s(x)|^2=\left(1-\frac{j_0}{j}\right)
\p(x)
+\frac{j}{j_0}h(x)-\frac{\log C_4}{2j}.
$$
We infer
\begin{equation}
\f_{j,j_0}(x) \geq \left(1-\frac{j_0}{j}\right)\f(x)
-\frac{\log C_4}{2j} \geq \f(x)-\frac{\log C_4}{2j}
\end{equation}
since $\f \leq 0$ on $X$. It follows from (3) and (4) that 
$\f_j \rightarrow \f$ in $L^1(X)$.

{\bf Step 2.}
We now show, following [18] that $(\f_{j,j_0})_j$ is almost subadditive.
Let $s \in \Gamma(X,L^{j_1+j_2})$ with
$$
\int_X |s|^2 e^{-2 h_{j_1+j_2,j_0}} dV_{\o} \leq 1.
$$
We may view $s$ as the restriction to the diagonal $\Delta$ of $X \times X$
of a section $S \in \Gamma(X \times X,L_1^{j_1} \otimes L_2^{j_2})$,
where $L_i=\pi_i^* L$ and $\pi:X \times X \rightarrow X$ denotes the projection
onto the $i^{th}$ factor, $i=1,2$.
Consider the Bergman spaces
\begin{eqnarray*}
\lefteqn{
{\mathcal H}_{j_1,j_2,j_0}:= \left\{ S \in 
\Gamma(X \times X,L_1^{j_1} \otimes L_2^{j_2}) \, / \, \right.} \\
&\hskip1.5cm& 
\left. \int_{X \times X} |S|^2 e^{-2h_{j_1,j_0/2}(x)-2h_{j_2,j_0/2}(y)}
dV_{\o_1}(x) dV_{\o_2}(y) <+\infty \right\},
\end{eqnarray*}
where $\o_i=\pi^* \o$. It follows from the Ohsawa-Takegoshi-Manivel 
$L^2$-extension theorem [14] that there exists 
$S \in \Gamma(X \times X,L_1^{j_1} \otimes L_2^{j_2})$ such that
$S_{|\Delta}=s$ and 
$$
\int_{X \times X} |S|^2 e^{-2h_{j_1,j_0/2}-2h_{j_2,j_0/2}}
dV_{\o_1} dV_{\o_2} \leq C_5
\int_X |s|^2 e^{-2 h_{j_1+j_2,j_0}} dV_{\o} \leq C_5,
$$
where $C_5$ only depends on the dimension $n=\dim_{\C} X$.
Observe that 
$\{\sigma_{l_1}^{(j_1,j_0/2)}(x) \cdot 
\sigma_{l_1}^{(j_1,j_0/2)}(y)\}_{l_1,l_2}$
forms an orthonormal basis of ${\mathcal H}_{j_1,j_2,j_0}$, thus
$$
S(x,y)=\sum_{l_1,l_2} c_{l_1,l_2} \sigma_{l_1}^{(j_1,j_0/2)}(x)
\sigma_{l_2}^{(j_2,j_0/2)}(y)
$$
with
$\sum |c_{l_1,l_2}|^2 \leq C_5$. It follows therefore from Cauchy-Schwarz
inequality that
$$
|s(x)|^2=|S(x,x)|^2 \leq C_5 \sum_{l_1} |\sigma_{l_1}^{(j_1,j_0/2)}(x)|^2
\sum_{l_2} |\sigma_{l_2}^{(j_2,j_0/2)}(y)|^2,
$$
which yields 
$$
\f_{j_1+j_2,j_0} \leq \frac{\log C_5}{2(j_1+j_2)}+
\frac{j_1}{j_1+j_2} \f_{j_1,j_0/2}+
\frac{j_2}{j_1+j_2} \f_{j_2,j_0/2}.
$$
Note finally that $\f_{j,j_0/2} \leq \f_{j,j_0}$ since $\f=\p-h \leq 0$,
therefore $\hat{\f}_j:=\f_{2^j,j_0}+2^{-j-2}\log C_5$ is decreasing.

{\bf Step3.} It remains to make $\hat{\f}_j$ smooth. Indeed it has all the
other required properties: it is decreasing and by Step 1 we have for all
$x \in X$,
\begin{equation}
\f(x) \leq \hat{\f}_j(x) \leq 
(1-j_02^{-j}) \sup_{B(x,r)} \f+\frac{C_6-n\log r}{2^j},
\end{equation}
so that $\hat{\f}_j \rightarrow \f$.
Let $\sigma_{1+s_{2^j}}^{(2^j)},\ldots,\sigma_{N_j}^{(2^j)} \in
\Gamma(X,L^{2^j})$ be such that $(\sigma_l^{(2^j)})_l$ is a basis
of $\Gamma(X,L^{2^j})$ and set
$$
\f_j:=\frac{1}{2^{j+1}} \log \left[ 
\sum_{l=1}^{s_{2^j}} |\sigma_l^{(2^j)}|^2
+\e_j \sum_{l=1+s_{2^j}}^{N_j} |\sigma_l^{(2^j)}|^2 \right]
+\frac{\log C_5}{2^{j+2}}-h.
$$
Clearly $\f_j \in PSH(X,\o)$. Moreover $\f_j \in {\mathcal C}^{\infty}(X)$
because $L^{2^j}$ is very ample if $j$ is large enough (hence we can find,
for every $x \in X$, a holomorphic section of $L^{2^j}$ on $X$
which does not vanish at $x$). Finally we can choose $\e_j>0$ that decrease
so fast to zero that $(\f_j)$ is still decreasing and converges to $\f$.
\end{proof}

\begin{cor}
Let $\o$ be a K\"ahler form on a {\bf projective} algebraic manifold X.
Then there exists $A \geq 1$ such that for every $\f \in PSH(X,\o)$,
we can find $\f_j \in PSH(X,A\o) \cap {\mathcal C}^{\infty}(X)$ which decrease
towards $\f$.
\end{cor}

\begin{proof}
Let $\f \in PSH(X,\o)$. Since $X$ is projective, we can find a Hodge
form $\o'$. Then $C^{-1} \o' \leq \o \leq C \o'$ for some constant $C \geq 1$.
Since $PSH(X,\o) \subset PSH(X,C\o')$, it follows from the previous theorem
that we can find $\f_j \in PSH(X,C \o') \cap {\mathcal C}^{\infty}(X)$
that decrease towards $\f$. Now the result follows from 
$PSH(X,C \o') \subset PSH(X,A \o)$with $A=C^2$.
\end{proof}

\begin{rqe}
When $X$ is merely K\"ahler, the above result still holds but the proof is
far more intricate.
We refer the reader to Demailly's papers for a proof.
\end{rqe}

\vskip .2cm

Vincent Guedj \& Ahmed Zeriahi

Laboratoire Emile Picard

UMR 5580, Universit\'e Paul Sabatier

118 route de Narbonne

31062 TOULOUSE Cedex 04 (FRANCE)

guedj@picard.ups-tlse.fr

zeriahi@picard.ups-tlse.fr

\end{document}